\documentclass[12pt]{article}

\usepackage{amsmath}
\usepackage{amssymb}
\usepackage{amsfonts}

\usepackage{amsthm}

\newtheorem{theo}{Theorem}

\newtheorem{prop}{Proposition}
\newtheorem{lemm}{Lemma}

\newcommand{\lbl}{\label}

\newcommand{\eq}[1]{$(\ref{#1})$}

\newcommand{\ER}{{Erd\H{o}s-R\'enyi~}}

\newcommand{\Po}{{\cal P}}
\newcommand{\Prv}{{Z}}
\newcommand{\PPrv}{{Z}}

\newcommand{\cM}{{\cal M}}

\def\N{\mathbb{N}}

\def\Real{\mathbb{R}}
\def\R{\mathbb{R}}
\def\E{\mathbb{E}}
\def\BX{\mathbb{X}}
\def\BM{\mathbb{M}}
\def\cX{\mathcal{X}}
\def\cS{\mathcal{S}}
\def\cY{\mathcal{Y}}
\def\cM{\mathcal{M}}
\def\bm{{\bf m}}
\def\bt{{\bf t}}
\def\bu{{\bf u}}
\def\cV{\mathcal{V}}
\def\cF{\mathcal{F}}

\def\Pr{\mathbb{P}}

\def\pmax{\overline{\phi}}

\def\0{{\bf 0}}
\def\bx{{\bf x}}
\def\by{{\bf y}}
\def\bz{{\bf z}}
\def\bX{{\bf X}}
\def\bu{{\bf u}}

\def\bS{{\bf S}}

\def\G{{\cal G}}
\def\De{D}
\def\Leb{{\rm Leb}}

\renewcommand{\E}{\mathbb E \,}

\newcommand{\tod}{\stackrel{{\cal D}}{\longrightarrow}}

\newcommand{\toP}{\stackrel{{P}}{\longrightarrow}}

\newcommand{\eqco}{\setcounter{equation}{0}}
\newcommand{\thco}{\setcounter{theo}{0}}
\newcommand{\prco}{\setcounter{prop}{0}}
\newcommand{\laco}{\setcounter{lemm}{0}}
\newcommand{\coco}{\setcounter{coro}{0}}
\newcommand{\cjco}{\setcounter{conj}{0}}

\newcommand{\deco}{\setcounter{defn}{0}}
\newcommand{\allco}{\eqco  \thco \prco \laco \coco \cjco \deco}
\newcommand{\LL}{{\cal L}}

\newcommand{\Y}{{\cal Y}}

\newcommand{\F}{{\cal F}}
\renewcommand{\G}{{\cal G}}

\newcommand{\tF}{{\tilde{F}}}
\newcommand{\tf}{{\tilde{f}}}
\newcommand{\tGG}{{\tilde{\G}}}

\newcommand{\phio}{\phi_0}

\newcommand{\eps}{\varepsilon}
\def\bdm{\begin{displaymath}}
\newcommand{\edm}{\end{displaymath}}
\def\benu{\begin{enumerate}}
\def\eenu{\end{enumerate}}
\def\beqn{\begin{equation}}
\def\eeqn{\end{equation}}
\def\be{\begin{equation}}
\def\ee{\end{equation}}
\def\bea{\begin{eqnarray}}
\def\eea{\end{eqnarray}}
\newcommand{\bean}{\begin{eqnarray*}}
\newcommand{\eean}{\end{eqnarray*}}
\newcommand{\bear}{\begin{eqnarray}}
\newcommand{\eear}{\end{eqnarray}}

\renewcommand{\epsilon}{\varepsilon}

\def\R{\mathbb{R}}

\def\qed{\hfill\hbox{${\vcenter{\vbox{
    \hrule height 0.4pt\hbox{\vrule width 0.4pt height 6pt
    \kern5pt\vrule width 0.4pt}\hrule height 0.4pt}}}$}}

\def\la{{\lambda}}

\begin{document}
\title{\bf 
Inhomogeneous random graphs,
isolated vertices, 
and Poisson approximation
}

\author{Mathew D. Penrose$^{1}$ 
\\
{\normalsize{\em University of Bath}} }

\maketitle


 \footnotetext{ $~^1$ Department of
Mathematical Sciences, University of Bath, Bath BA2 7AY, United
Kingdom: {\texttt m.d.penrose@bath.ac.uk} }
%




\begin{abstract}
Consider a graph on randomly scattered points
in an arbitrary space,
with any two points $x,y$ connected with
probability $\phi(x,y)$. Suppose the
number of points is large but the
mean number of isolated points is $O(1)$.
We give general criteria for the latter
to be approximately Poisson distributed. 
More generally, we consider the number of
vertices of fixed degree, the number
of components of fixed order, and the number of edges.
We use a general result on Poisson approximation
by Stein's method for a set of points selected
from a Poisson point process. This method also gives
a good Poisson approximation for U-statistics of a
Poisson process.
\end{abstract}



\section{Introduction}
\label{sec.intro}
In the {\em inhomogeneous random graph} (IRG),
each vertex  has one of several 
 possible types or {\em states}, where the
space of possible states may be infinite.
Given the states of the vertices, each possible edge
is present with a probability that depends on the states of
the two endpoints, independently of the other edges.
This provides a very flexible class of random graph models;
for example, the state of  a vertex could represent
spatial location, or it could represent  the time
at which a vertex is born, for a randomly evolving 
graph.

Such models are popular 
 in statistical network modelling,  where they go under names 
such as {\em stochastic block model} \cite{Blockref}
(in which case the state space is usually taken to be finite)
and {\em latent space model} 
\cite{Hoff}
or {\em latent variable model} \cite{Latentvarref}.
 The IRG terminology dates from
  \cite{Soderberg}, and is prevalent in the more probabilistic
literature, where the model has been studied in
depth, for example in \cite{BJR} and 
\cite{BHL}.
Much of this literature is concerned
with the birth of the giant component, but full connectivity 
has also been studied in \cite{DF}.

 For our purposes, the IRG is defined as follows.
Let $(\BX,\cF,\mu)$ be a probability space (the state space).
For $s >0$ let $\Po_s$ be the (random) set of points of a Poisson
point process on $\BX$ with mean measure (i.e., intensity
measure) $s \mu$.
Also, for $n \in \N$, let $\cX_n$ be the
binomial point process consisting of $n$ independent random elements
of $\BX$ which common distribution $\mu$.

Suppose $\phi: \BX \times \BX \to [0,1]$  
is a measurable symmetric function; 
we call such a function
a {\em connection function}.
Given finite $\cV \subset \BX$ (possibly with multiplicity), let
$G(\cV,\phi)$ be 
 the random (undirected) graph
with vertex set $\cV$, with
 each pair $\{x,y\}$ of points of $\cV$
 connected by an edge with probability 
$\phi(x,y)$, independently of all other pairs.
We are particularly interested in $G(\Po_s,\phi)$ and
$G(\cX_n,\phi)$; 
we define these graphs more formally
in Section \ref{secconstruct}.

In the special case where $\BX$ is a region
of Euclidean space and $\phi(x,y)$ is
determined by the displacement $x-y$ (typically via $\|x-y\|$,
where $\|\cdot \|$ denotes the Euclidean norm),
the IRG is also known as the 
 {\em soft random geometric graph} \cite{Soft} or 
{\em random connection model} (RCM) \cite{MRBk}.
If, in the Euclidean setting, we have 
$\phi(x,y)= {\bf 1}_{\{\|x-y\| \leq r \}}$
the IRG is known as the {\em random geometric graph} (RGG) 
\cite{PBk} or
{\em Gilbert graph.}
These models are important in applications to wireless communications;
see for example \cite{DPS,FM,GK,MaoAnderson,YWLH}.

For any graph $G$ and any $j \in \N_0:= \N \cup \{0\}$, let $\De_j(G)$
denote the number of vertices in $G$ of degree $j$; also
 set $\De_{\leq j}(G) := \sum_{i=0}^j \De_i(G)$.
In particular, $\De_0(G)$ is the number of isolated vertices. 
Of interest is the question of whether  $\De_{\leq k}(G(\Po_s,\phi))$ 
and  $\De_{\leq k}(G(\cX_{\lfloor s \rfloor },\phi))$ are
 approximately Poisson distributed for
$s$ large, with $k$ fixed, e.g. $k=0$. 
One reason for interest is that $\De_{0}(G) =0$
(respectively $\De_{\leq k}(G) =0$) is clearly a necessary
condition for $G$ to be connected (resp. $(k+1)$-connected,
assuming $G$ has at least $k+2$ vertices), 
and for  many choices of $\BX$ and $\phi_n$
 this condition is asymptotically sufficient 
(in probability) when $G=G(\cX_n,\phi_n) $ with $n$ large; see
 \cite{Iyer,Hsing,Soft,kcon,Norm}.
 In such cases,
$\Pr [ \De_{\leq k}(G(\cX_n,\phi_n)) =0]$ is
a good approximation for the probability that
$G(\cX_n,\phi_n)$ is $(k+1)$-connected,
so if we can estimate the former via Poisson
approximation, we may have a useful estimate for
the probability of $(k+1)$-connectivity.

The references just mentioned show that in many cases
where $s$ is large and $\phi $ is chosen so that
$\E \De_j(G(\Po_s,\phi)) = O(1)$, the distributions of 
$\De_j(G(\Po_s,\phi)) $ and of
$\De_j(G(\cX_{\lfloor s \rfloor},\phi)) $
 are approximately 
Poisson. In this paper we give a general criterion for this
to be the case, without making any geometrical or topological assumptions on
the space $\BX$ (Theorem \ref{Th.k}). We give similar results for the number
of components of order $k$ with $k$ fixed (Theorems \ref{Th.moments}
and \ref{Th.cpts}),
and for the number of edges (Theorem \ref{thclump}).
These theorems also incorporate asymptotic normality, when
the mean of the variable in question grows slowly as a function of $s$.

We prove Theorems \ref{Th.k}, \ref{Th.moments}, \ref{Th.cpts}
 and \ref{thclump} in
Sections \ref{secproveThk}, \ref{secGamma} and \ref{secGn}.
The proofs use
Theorem \ref{thSteink}, a general result on Poisson approximation for
functionals of Poisson processes, which is of independent interest.
 This theorem also gives us further results on the Poisson approximation 
for the number of edges under different assumptions
from those of Theorem \ref{thclump}. See Theorem \ref{th.edge}
and the subsequent discussion.
It also provides an alternative way to derive
(and slightly improve) a result of \cite{DST}  on
Poisson approximation of $U$-statistics of a Poisson process; see
Section \ref{secU}.

\section{Results on inhomogeneous random graphs}
\label{secresults1}
\allco


Let $\Phi$ be the class of all measurable symmetric functions 
from $\BX \times \BX $ to $ [0,1]$,
and for $\eps >0$ let $\Phi_\eps$ be the class of 
 all $\phi \in \Phi $ that satisfy
$$
\inf_{x \in \BX } \left( \int_\BX \phi(x,y) \mu( dy) \right) 
\geq \eps
\sup_{x \in \BX} \left( \int_\BX  \phi(x,y) \mu(dy)   \right).
$$
If $\phi \in \Phi_\eps$ we say the connection function
$\phi$ is $\eps$-{\em homogeneous}.
Note that $1$-homogeneity, according to our definition,
is the same as homogeneity as defined in \cite{DF}.

Several of our results require $\eps$-homogeneity. 
There are many interesting classes of connection
function which satisfy this condition. 
 For example, suppose $\BX$ is a bounded convex region
in Euclidean space $\R^d$, and $\mu$ has a density (with
respect to Lebesgue measure on that region) that is bounded away from 0 and
infinity. Then there exists $\eps >0$ such that
 all connection functions of the form $\phi(x,y) = \psi(|x-y|)$ 
with $\psi $ nonincreasing lie in $\Phi_\eps$. In
 particular, we do not require any exponential
decay condition on $\psi$, such as that imposed in
\cite{Soft}.

  For $\phi \in \Phi$, set $\pmax = \sup_{x,y} \phi (x,y)$.  
As well as $\eps$-homogeneity, some of our results also require
that $\pmax$ be bounded away from 1. This condition is annoying since
it excludes from consideration the standard RGG,
and also some cases of the IRG which have been considered
elsewhere in the literature \cite{BJR,DF}.  On the other hand,
 `soft' versions of these models, which do satisfy this condition,
are easily defined and  arguably  
 will often be reasonable from the point of view of applications
(we shall discuss this further below).  Without this condition,
it may be unrealistic to expect to find a simple
argument for Poisson approximation
of $D_0$ covering a large class of RGG densities 
without getting involved in geometrical details of any particular
probability  density or even assuming any Euclidean structure
at all, which is what our condition
on $\pmax$ allows us to do.

Suppose we have connection functions  $\phi_s$ defined for each $s >0$
with $\E \De_{\leq j}(G(\Po_s,\phi_s)) $ bounded as $s \to \infty$,
for some fixed $j \in \N_0$.
See (\ref{ENk}) below for a formula for  $\E \De_{j}(G(\Po_s,\phi_s)) $.
Our first result shows that $\De_j(G(\Po_s,\phi_s))$ 
 and $\De_j(G(\cX_{\lfloor s \rfloor } ,\phi_{\lfloor s \rfloor } ))$ 
are approximately Poisson 
and $\De_k(G(\Po_s,\phi_s))$ is approximately normal for $k > j$, 
for large $s$ under the condition that
 $\phi_s \in \Phi_\eps$ for all $s$ and $\pmax_s $ is bounded away from 1.  

For $\alpha \in (0,\infty)$,
let ${\Prv}_\alpha$ denote a random variable having
 the Poisson distribution with parameter $\alpha$.
Let ${\cal N}$ denote a random variable having
 the standard normal distribution in $\R$,
i.e. with probability density function 
$(2\pi)^{-1/2} \exp(- x^2/2)$, $x \in \R$.
For any graph $G$ we define $D_{\leq -1}(G):=0$.
For any two $(0,\infty)$-valued functions $u(s)$ and $v(s)$,
we say $u (s) = O (v(s))$ if $\limsup_{s \to \infty} u(s)/v(s) <\infty$,
and $u(s) =o(v(s))$ if $\limsup_{s \to \infty} u(s)/v(s) =0$, and
 $u(s) = \Theta(v(s))$ if
 $u(s)= O(v(s))$ and $v(s) = O(u(s))$.


\begin{theo}
\label{Th.k}
Let $j,k \in \N_0$, $\eps >0$
and $\phi_s \in \Phi_\eps$ with $\overline{\phi}_s \leq 1- \eps$
for $s \in (0,\infty)$, satisfying
$\lim_{s \to \infty}
 \E \De_{\leq j}(G(\Po_s,\phi_s)) = \alpha \in (0,\infty)$. Then
as $s \to \infty$,

\bea
\De_j(G(\Po_s,\phi_s)) \tod {\Prv}_\alpha;  ~~~~~
 \De_j(G(\cX_{\lfloor s \rfloor},\phi_{\lfloor s \rfloor})) 
\tod {\Prv}_\alpha,
\label{eqpolim}
\eea
and also
\bea
\De_{\leq j-1}(G(\Po_s,\phi_s)) \toP 0; ~~~~~ 
 \De_{\leq j-1}(G(\cX_{\lfloor s \rfloor},\phi_{\lfloor s \rfloor})) \toP 0,
\label{eq0lim}
\eea
and moreover
 \bea
\E[ \De_k(G(\Po_s,\phi_s)) ] = \Theta((\log s)^{k-j})
\label{EDiDj}
\eea
and if also $k >j$ then 
\bea
 \frac{ \De_k(G(\Po_s,\phi_s)) - \E[ \De_k(G(\Po_s,\phi_s)) ] 
}{ 
 \sqrt{ \E[ \De_k(G(\Po_s,\phi_s)) ] } } \tod {\cal N}.
\label{DiN}
\eea
\end{theo}
It is interesting to compare the conclusions
(\ref{eqpolim}) and (\ref{eq0lim}) of
 this result with the example on page 55 of
\cite{PBk}. In that case, for a certain sequence of  RGGs
one can arrange for the number of vertices of degree
2 to be asymptotically {\em compound} Poisson, whereas
here it is asymptotically Poisson.

For $k \in \N$, and for any graph $G$, we refer to the components
 of $G$ of order $k$ (i.e., with $k$ vertices)
 as the {\em $k$-components} of $G$.
 Let $N_k(G)$ denote the number of $k$-components 
in $G$.  In particular, $N_1(G)= D_0(G)$.

Suppose we have connection functions  $\phi_s$ defined for each $s >0$
with $\E N_k(G(\Po_s,\phi_s)) $ bounded as $s \to \infty$, for
some fixed $k \in \N$. 
Theorem \ref{Th.cpts} below shows that $N_k(G(\Po_s,\phi_s))$
and  $N_k(G(\cX_{\lfloor s \rfloor},\phi_{\lfloor s \rfloor}))$ are
 approximately Poisson under
the condition that  $\phi_s \in \Phi_\eps$
and $\pmax_s \leq 1-\eps$
for all $s$ and some fixed $\eps >0$; even without this condition, 
Theorem \ref{Th.moments} shows that
the Poisson approximation holds provided
that $\pmax_s = o(1/\log s)$ 
 (in the case of $G(\Po_s,\phi_s)$), or provided
that $\pmax_{\lfloor s \rfloor} = o(s^{-1/2})$ 
 (in the case of $G(\cX_{\lfloor s \rfloor},\phi_{\lfloor s \rfloor})$).
These rates of decay imposed  on $\pmax_s$ are
 significantly milder than the condition $p_n = \Theta((\log n)/n)$
for the \ER random graph $G(n,p_n)$ to be 
at the threshold for having no isolated vertices.

We also give a result on asymptotic normality, when
 $\E N_{k}(G(\Po_s,\phi_s)) $ grows slowly as $s \to \infty$.
All asymptotics in the next two theorems are as $s \to \infty$.

\begin{theo}
\label{Th.moments}
Let $k \in \N$ and
 $\phi_s \in \Phi$ for $s > 0$, with
$\E N_{k}(G(\Po_s,\phi_s)) \to \alpha \in (0,\infty)$.

(a) If 
 $\overline{\phi}_s = o(1/\log s)$, then
$N_k(G(\Po_s,\phi_s)) \tod {\Prv}_\alpha$.

(b) If
 $\overline{\phi}_s = o(s^{-1/2})$,
then
$N_k(G(\cX_{\lfloor s\rfloor},\phi_{\lfloor s\rfloor})) 
\tod {\Prv}_\alpha$.
\end{theo}

\begin{theo}
\label{Th.cpts}
Let $k \in \N, \eps >0$ and
 $\phi_s \in \Phi_\eps$ for $s > 0$, with
$\pmax_s \leq 1-\eps$
for all $s$. 

(a) 
If $ \E N_{k}(G(\Po_s,\phi_s)) \to \alpha \in (0,\infty)$,
 then 
\bea
N_k (G(\Po_s,\phi_s)) \tod  {\Prv}_\alpha; ~~~~~ 
N_k (G(\cX_{\lfloor s \rfloor},\phi_{\lfloor s \rfloor})) \tod  
{\Prv}_\alpha. 
\label{limGPo}
\label{limGBi}
\eea

(b)
If 
  $\E N_{k}(G(\Po_s,\phi_s)) \to \infty$,
but $\E N_{k}(G(\Po_s,\phi_s)) = o(s)$,
then  setting $\alpha_s = \E N_{k}(G(\Po_s,\phi_s)) $,
we have that $ (N_{k}(G(\Po_s,\phi_s)) - \alpha_s) /\sqrt{\alpha_s}
 \tod {\cal N}$.
\end{theo}

Among other things, the case $k=1$ of
Theorem \ref{Th.moments} (a) 
generalizes Lemma 3.2 of \cite{Soft} to 
a more general class  of $(\BX, \phi_s)$
than is considered in \cite{Soft}.

Given a connection function $\phi \in \Phi$, let us define
\bea
\kappa(\phi) := \sup_{x \in \BX} \int_{\BX} \phi(x,y) \mu(dy).
\label{kappadef}
\eea
If $\phi$ is $\eps$-homogeneous for fixed $\eps >0$, then for any $x \in \BX$
the expected degree of a vertex of $G(\Po_s,\phi)$ located at $x$
is of the order of  $s \kappa(\phi)$.
Our proof of 
Theorems \ref{Th.k} and \ref{Th.cpts} actually shows that
under conditions of $\eps$-homogeneity and $\overline{\phi}_s \leq 1-\eps$,
if $s\kappa(\phi_s) \to \infty$ (so the `typical degrees' become large)
then for any fixed $j \geq 0$ and $k \geq 2$,
both $D_j$ and $N_k$ are approximately Poisson distributed (so
if the mean of this Poisson distribution goes to infinity then
they are approximately normally distributed).
Moreover, under the conditions of Theorem \ref{Th.k} 
or Theorem \ref{Th.cpts} (a),
we shall show that 
$s\kappa(\phi_s) = \Theta(\log s)$ (see (\ref{Theta2cond})
and (\ref{Thetacond})).

 In Theorem \ref{Th.k} we do
 not address asymptotic normality  $D_j(G(\Po_s,\phi_s))$ when
its expected value grows to infinity more slowly than $s$,
except for the case $j=0$ which is covered by the case $k=1$
of Theorem \ref{Th.cpts} (b).
In attempting to adapt the proof of 
the latter to the case $j>0$, the difficulty is that 
 in general the expression (\ref{EtU})
does not vanish in the case where $s \kappa (\phi_s)$ tends to zero.
We conjecture that asymptotic normality  of $D_j(G(\Po_s,\phi_s))$ can
be proved by other means but this is beyond the
scope of this paper.

Our next result is concerned with Poisson or normal approximation for
(a generalization of)
the number of edges of $G(\Po_s,\phi_s)$. 
 This is of interest in itself, and will also be of use in the proof of
Theorem \ref{Th.cpts}.
Given $k \in \N$ with $k \geq 2$, for any graph $G$
let $H_k(G)$ denote the
number of connected induced subgraphs  
of $G$ of order $k$ (so in particular, $H_2(G)$ is
the number of edges.)
Again, asymptotics are as $s \to \infty$.
\begin{theo}
\label{thclump}
Let $\eps >0$,
 $\phi_s \in \Phi_\eps$ for all $s >0$,
 and $k \in \N$ with $k \geq 2$.

(a) 
If $\E H_k(G(\Po_s,\phi_s)) \to \alpha \in (0,\infty)$,
then $ H_k(G(\Po_s,\phi_s)) \tod {\Prv}_\alpha $.

(b)
If $\E H_k(G(\Po_s,\phi_s)) \to \infty$ but 
$\E H_k(G(\Po_s,\phi_s)) = o(s)$, then setting $\alpha_s :=
\E H_k(G(\Po_s,\phi_s))$, we have
$(H_k(G(\Po_s,\phi_s)) - \alpha_s)/\sqrt{\alpha_s } \tod {\cal N} $.
\end{theo}

The proof of Theorems \ref{Th.k}, \ref{Th.cpts} 
 and \ref{thclump}
provides information about the rates of convergence in these results.
Theorem \ref{Th.moments}  will be proved by the method of moments,
which does not provide any information about rates.
The reason part (b) of that result requires 
a stronger condition on $\overline{\phi}_s$ than  part (a) does,
is because when bounding these moments from below, for part 
(a) we have a sequence of the form $\exp(-n x_n)$ (with  
 $x_n$ bounded by a constant times $\overline{\phi}_n$), while
for (b) we have a sequence of the form $(1-x_n)^n$ 
which  is asymptotic  to $\exp(-nx_n)$ only when $x_n = o(n^{-1/2})$

We now discuss some of the literature related to these results.
Much of this concerns the RGG.
Suppose $(\BX,\cF,\mu)$ is the $d$-dimensional Euclidean space 
 with the Borel $\sigma$-algebra and with
$\mu$ having a density $f$ with respect to Lebesgue measure, while
the connection function is given by $\phi_s(x,y) =
 {\bf 1}_{\{|x-y| \leq r_s\}}$, 
with $r_s$ chosen in such a way that
$\E D_0(G(\Po_s,\phi_s)) \to \alpha \in (0,\infty)$.
It has been shown that $D_0 \tod {\Prv}_\alpha$ when $\mu$ is
uniform on the unit cube \cite{kcon}, and for 
certain special types of density with unbounded support \cite{Iyer,Hsing,Norm}.
However, $D_0$ is not always asymptotically Poisson;
see the last paragraph of \cite[Section 2]{Hsing}, where it
is suggested that the Poisson limit 
is `the exception rather than the rule' in dimension $d=1$.
We are a long way from having any complete characterisation
of distributions for which the number of isolated vertices in
the RGG is asymptotically Poisson.


 
Suppose $(\BX,\cF,\mu)$ is the $d$-dimensional unit cube
equipped with Lebesgue measure.
If $\phi_s$ is of the form $\phi_s(x,y) =
 p_s {\bf 1}_{\{|x-y| \leq r_s\}}$, with $p_s$ and $r_s$
 chosen so that $\E D_0(G(\Po_s,\phi_s)) \to 
\alpha$, then $D_0 \tod {\Prv}_\alpha$.
More generally, if $\phi_s$ is of the form 
$\phi_s(x,y) = \psi_s(|x-y|)$ with $\psi_s$  
a decreasing function satisfying an exponential
decay condition, then the Poisson convergence of $D_0$ 
is known to hold; see \cite{Soft}.
Our results enable us to relax the exponential decay condition
and allow for other distributions, for example with
 density bounded away from zero and infinity  on a convex compact region
in $\R^d$.

Here is a simple example of a sequence of connection functions $\phi_s$
where 
\linebreak
 $\De_0(G(\Po_s,\phi_s)) $ has bounded mean, but is
{\em not} 
asymptotically Poisson.  Let $(\BX,\F,\mu)$ be
 the unit interval equipped with Lebesgue measure,
and set $\phi_s (x,y) =1$  if $\max(x,y) \leq s^{-1}$
or $\min(x,y) > s^{-1}$,  otherwise setting
$\phi_s(x,y) =0$. Then in the large-$s$  limit the random
variable $\De_0(G(\Po_s,\phi_s))$
is Bernoulli distributed with parameter $e^{-1}$, not Poisson distributed. 
The condition of  $\eps$-homogeneity, which appears in
many of our results, rules out this sort of example.

 Devroye and Fraiman \cite{DF} consider $D_0(G(\cX_n,\phi_n))$  on
a general space $\BX$ in the case where $\phi_n = (a(n) \phio) \wedge 1$
for a fixed function $\phio: \BX \times \BX \to \R_+$
(such a function $\phio$ is called a {\em kernel}), and some sequence
$a(n)$ (where $\wedge$ denotes minimum). This is a common assumption for
the IRG; see for example \cite{BJR}. 
Their results may be interpreted as saying that
(under certain conditions) the threshold
value of $a$ above which $G(\cX_n, a \phio \wedge 1)$ 
is free of isolated vertices (which may be viewed
as a random variable) satisfies a weak law of large numbers;
it is asymptotic to a constant times  $(\log n)/n$, in
probability.  They also derive a similar law of large numbers
 for the threshold for the graph to be connected.

When $\phio$ is bounded, our Theorem \ref{Th.moments} 
 (b) gives a possibility of
deriving convergence in distribution for this threshold
(suitably transformed).
For example, suppose $\BX =[0,1)$ and
$\mu$ has a density $f$ with
respect to  Lebesgue measure, with $f$ bounded away from zero. 
Suppose also that $\phio(x,y)$ is determined via the absolute
value of $x-y$ (mod 1) (so we are in the one-dimensional torus).
If $f \equiv 1$, or if $f$ is smooth with $ g(x) := \int \phio(x,y) f(y) (dy)$
having a unique minimum, then it should be  possible to
derive a distributional limit law for this threshold,
since in these cases it should be possible, for any $\alpha \in (0,\infty)$,
to determine a sequence $a_s$ such that 
 $\E D_0(G(\Po_s,  a_s \phio)) \to \alpha$, and then apply Theorem
\ref{Th.moments} (b).

If $\phio$ is unbounded but shift-invariant on the torus,
for example if $\phio(x,y) = ((x-y){\rm mod} ~1)^{-\gamma}$
for some fixed positive $\gamma$ and  $f \equiv 1$, 
then for any $\eps \in (0,1)$, it may  be possible
 to find a limiting distribution
for a suitable  transformation of the threshold value of $a$ above which
  $G(\cX_n,a \phio \wedge (1-\eps))$ is
  free of isolated vertices, now using Theorem \ref{Th.cpts} (a).
 Indeed, in this case the connection functions
are all $1$-homogeneous, and by using connection functions 
  $\phi_n = a_n \phio \wedge (1-\eps)$ rather than the more
standard $\phi_n = a_n \phio \wedge 1$,
we ensure that the extra condition
$\pmax_n \leq 1-\eps$ is also satisfied.

It would be interesting to fully work out and extend 
these examples; to get a similar results for the 
connectivity threshold;  and to improve the weak law of
\cite{DF} to a strong law.

In the present work we consider only undirected graphs.
Analogous directed graph  models can be defined similarly;
these have
been considered in the random geometric graphs
literature \cite{IT} and in the statistical
literature \cite{Hoff}.
 In this case, the connection function $\phi(x,y)$
represents the probability that there is a directed edge from
a vertex at $x$ to a vertex at $y$. It
 is not required to be symmetric, and can be adapted so that
all vertices have the same expected outdegree, as in \cite{IT}.
It would be interesting to try to
 derive a similar result to 
Theorem \ref{Th.k}  for the number of vertices with
outdegree $j$ in such a model.

\section{A general result on Poisson approximation}
\label{secstein}
\allco

Let $(\BM,\cM,\bm)$ be a probability space
(known as a {\em mark space}).
Assume that the probability measure $\bm$ on $\BM$ 
is {\em diffuse}, by which we mean that
there is a product measurable set $A \subset \BM \times \BM$
with $(\bm \otimes \bm)(A) =0$,
 such that the diagonal $\{(t,t):t \in \BM \}$
is contained in $A$. For example, if $\{t\} \in \cM$
and $\bm (\{t\} ) =0$ for all $t \in \BM$,
then $\bm$ is diffuse. 

Suppose on a suitable probability space that
we have a sequence $((X_i,T_i),i=1,2,3,\ldots)$ of
 independent identically
distributed random elements of $\BX \times \BM$ with common distribution
$\mu \otimes \bm$, and an independent unit rate 
Poisson counting process $(\PPrv_s,s > 0)$, so that the  random variable
 $\PPrv_s$
 has Poisson($s$) distribution for each $s$, and also a further
independent sequence $(\tau,\tau_1,\tau_2,\tau_3,\ldots)$ of independent  
 random elements of $\BM$ with common distribution $\bm$.
By our assumption that $\bm$ is diffuse,
 the values of $\tau,T_1,\tau_1,T_2,\tau_2,T_3,\tau_3,\ldots$ are almost
surely distinct.

A {\em finite point process} in $\BX$ is defined as
a random element of the space $\bS(\BX)$ of all finite
 subsets of $\BX$, where $\bS(\BX)$ is
 equipped with the smallest $\sigma$-algebra
$\cS(\BX)$ containing the sets $\{\xi \in\bS(\BX):\xi(B)=k\}$
for all  $B \in \cF $
 and all $k\in\N_0$, where
$\xi(B):= |\xi \cap B|$ and  $|\xi|$ denotes the
 number of elements of $\xi $.

A  finite point process in $\BX \times \BM$
is defined similarly  as a random element of the space
 $\bS:= \bS(\BX \times \BM)$, where $\BX \times \BM$ is
equipped with the product $\sigma$-algebra $\cF \otimes \cM$.
For $k \in \N$, let $\bS_k := \{ \xi \in \bS: |\xi| =k\}$.

Given $s >0$, $n \in \N$,  define the following point processes 
in $\BX \times \BM$:
\bea
\eta_s := \cup_{i=1}^{\PPrv_s} \{(X_i,T_i)\} ; ~~~~
\xi_n := \cup_{i=1}^n \{(X_i,T_i)\}.
\label{etaxi}
\eea
Then $\eta_s$ is 
a Poisson point process in $\BX \times \BM$
with mean measure $s \mu \times \bm$.
Similarly, $\xi_n$ is a  binomial point process in $\BX \times \BM$. 

Let $d_{TV}$ and $d_W$ denote total variation distance
and Wasserstein distance, respectively, between probability
measures on the nonnegative integers.
That is, for $\N_0$-valued random variables
$X,Y$ with distribution $\LL(X),\LL(Y)$ respectively,
we set 
\bean
d_{TV}(\LL(X),\LL(Y)) = \sup_{A \subset \N_0} 
(\Pr[X \in A] - \Pr[Y \in A]);
\\
d_W(\LL(X),\LL(Y)) = \sup\{|\E h(X) - \E h(Y) | :
\| \Delta h\|_\infty \leq 1\},
\eean
where for $h: \N_0 \to \R$ we set 
$\Delta h(i) := h(i+1) - h(i)$ for $i \in \N_0$,
and $\|h\|_\infty := \sup_{i \in \N_0} |h(i)|$

The following theorem is related to
a well-known result on the Poisson approximation
of a sum of Bernoulli random variables by
Stein's method via coupling
(Theorem II.24.3 of \cite{Lindvall}, or
 Theorem 1.B of \cite{BHJ}).  Here the terms in
the sum are themselves indexed by $k$-subsets of
the set of points of a (marked)
 Poisson point process.

Let $k \in \N$ and
let $f: \bS_k \times \bS \to \{0,1\}$ be a measurable function.
For $\xi \in \bS$, set
\bea
F(\xi) := 
\sum_{\{\psi \subset \xi:|\psi| =k\}}
 f \left(\psi,\xi \setminus \psi \right) . 
\label{Wdef}
\eea 
We can think of $f$ as a mechanism for selecting some
of the $k$-subsets
of $\xi$, and $F(\xi)$
as the total number of $k$-subsets  selected.

\begin{theo}
\label{thSteink}
Let $s >0$.
Let $W := F(\eta_s)$ with $\eta_s$ and $F$ as described above.
For $x_1,\ldots,x_k \in \BX$
set $p(x_1,\ldots,x_k) := \E f(\{(x_1,\tau_1),\ldots,(x_k,\tau_k)\}, \eta_s),$
and set $\lambda := s \mu$.

Suppose that $w: \BX^k \to [0,\infty)$ is a measurable function,
 and that for $\lambda^k$-almost every $\bx = (x_1,\ldots,x_k) \in \BX^k$
 with $p(x_1,\ldots,x_k) >0$
  we can find coupled random 
variables $U_{\bx},V_{\bx}$ such that
\begin{itemize}
\item
$
\LL(U_{\bx}) = \LL(W);
$
\item
$
\LL(1+ V_{\bx}) = \LL(F (
 \cup_{i=1}^k 
\{(x_i,\tau_i)\}
\cup
\eta_s )
 | f(\cup_{i=1}^k \{(x_i,\tau_i) \} ,\eta_s
) =1 ).
$
\item
$\E [ |U_{\bx} - V_{\bx}|] \leq w(\bx)$.
\end{itemize}
 Set $\alpha = \E W = (1/k!)\int p(\bx) \la^k ( d\bx)$.
 Then
\bea
d_{TV}(\LL(W),\LL({\Prv}_\alpha)) \leq \frac{ ( 1 \wedge \alpha^{-1}) }{k!} 
\int 
w(\bx)
 p(\bx) \la^k (d\bx),
\label{SteinkTV}
\eea
and
\bea
d_{W}(\LL(W),\LL({\Prv}_\alpha))
 \leq \frac{ 3 ( 1 \wedge \alpha^{-1/2})  }{k!}
\int 
w(\bx)
 p(\bx) \la^k (d\bx).
\label{SteinkW}
\eea
\end{theo}
The proof uses the
 (multivariate) {\em Mecke formula} 
(see e.g. \cite{LP}, or  \cite[Theorem 1.6]{PBk}),
which says that if
 $ g: \bS_k \times \bS \to \R$
is a bounded measurable function,
then
\bea
\E \sum_{\{\psi \subset \eta_s: |\psi | =k \}}
 g\left(  \psi, \eta_s \setminus \psi \right)
= \frac{1}{k!} \int \E \left[ g(\{(x_1,\tau_1),\ldots,(x_k,\tau_k) \} ,\eta_s
) \right]  
 \lambda^k(d\bx).
\label{eqMecke}
\eea
This fact gives us the assertion in the statement of
the theorem that $\E W = (1/k!) \int p(\bx) \lambda^k(d\bx)$.

Theorem \ref{thSteink} still holds in the case where the
measure $\lambda$ is taken to be  $\sigma$-finite but infinite,
and $\eta_s$ is replaced by a Poisson point process $\eta$ with
mean measure $\lambda \otimes \bm$.
The proof is essentially unchanged. \\

{\em Proof of Theorem \ref{thSteink}.}
Let $h: \N_0 \to \R$ be bounded.
Then 
\bean
\E[W h(W)] 
= \E \sum_{\{\psi \subset \eta_s: | \psi | =k\} } 
 f( \psi, \eta_s \setminus \psi)
h(F(\eta_s))
\\
= \frac{1}{k!} \int \E \left[ f(
\{(x_1,\tau_1),\ldots,(x_k,\tau_k) \} ,
\eta_s
)  
h\left(F
 \left(
 \cup_{i=1}^k \{(x_i,\tau_i)\}
 \cup
\eta_s
\right)
\right) \right] 
\lambda^k(d\bx)
\\
= \frac{1}{k!} \int 
 \E\left[  \left. h\left(F\left(
 \cup_{i=1}^k 
\{(x_i,\tau_i)\} 
 \cup
\eta_s \right)
\right) 
 \right| f \left(\cup_{i=1}^k\{(x_i,\tau_i)\},\eta_s
\right) =1 \right] 
p(\bx)
\lambda^k (d\bx).
\eean
Also, 
$\E [\alpha h(W+1) ] = 
(1/k!) \int \E h(W+1) p(\bx) \la^k(d\bx)$,
 and therefore
\bean
|\E [\alpha h(W+1) -  W h(W)] |  
\leq \frac{1}{k!} \int
 p(\bx) \lambda^k(d\bx)
~~~~~~~~~~~~~~~~~~~~~~~~~~~~~
~~~~~~~~~~
\\
\times
 \left| \E h(W+1)  - \E \left[ 
  \left. h\left(F\left(
 \cup_{i=1}^k \{(x_i,\tau_i)\}
\cup
\eta_s 
\right)
 \right) \right| 
f\left(\cup_{i=1}^k\{(x_i,\tau_i)\},\eta_s
\right) =1 \right] 
\right|
.
\eean
Since $|h(i)-h(j)| \leq \|\Delta h\|_\infty \cdot |i-j|$ 
for $i,j \in \N_0$,
 we obtain for each $\bx = (x_1,\ldots,x_k) \in \BX^k$ that
\bean
\left| \E h(W+1)-\E
    \left[ \left. h \left( F
           \left( \cup_{i=1}^k  \{(x_i,\tau_i)\} \cup  \eta_s \right) \right)  
          \right|
 f \left( \cup_{i=1}^k \{(x_i,\tau_i)\},\eta_s \right) =1  \right] \right|
\\
\leq
|\E h(U_\bx+1)  - \E  h( V_\bx+1 ) | 
\\
\leq \|\Delta h \|_\infty \E |U_\bx - V_\bx | 
\leq \| \Delta h \|_\infty w(\bx), 
\eean
and therefore
\bean
| \E [ \alpha h(W+1) - Wh(W) ]| 
\leq \frac{\| \Delta h\|_\infty}{k!} \int w (\bx )  p(\bx) \lambda^k(d\bx).
\eean
Given $A \subset \N_0$, set
$g= {\bf 1}_A$ and choose $h:\N_0 \to \Real$ so that $h(0)=0$ and
\bea
\alpha h(i+1) - i h(i) = g (i) - \E[G({\Prv}_\alpha ) ] , i \in \N_0.
\label{hfromg2}
\eea
Then (see Lemma 1.1.1 of \cite{BHJ}) $h$ is bounded and
$
\| \Delta h\|_\infty \leq 1 \wedge \alpha^{-1},
$
and hence
$$
|\Pr[W \in A] - \Pr[ {\Prv}_\alpha \in A] | 
\leq \frac{ (1 \wedge \alpha^{-1})}{k!}
\int 
w(\bx) p(\bx) \la^k(d\bx).
$$
 The result (\ref{SteinkTV}) follows.

One obtains (\ref{SteinkW}) similarly by 
choosing, for any given $g$ with $\|\Delta g \|_\infty \leq 1$,
 a solution $h$ to (\ref{hfromg2})
with $h(0)=0$, and using Lemma 1.1.5 of \cite{BHJ}.
 $\qed$ \\

We now give an overview of how we shall use Theorem 3.1 to prove
Theorems \ref{Th.k} and \ref{Th.cpts}.
As explained in the next section, we may view $G(\Po_s,\phi_s)$ as
being determined by a  marked Poisson point process in $\BX$, i.e.
a Poisson point process in a product space $\bX \times \BM^*$ (where
$\BM^*$ is a certain mark space). 

Then the functional $D_k (G(\Po_s,\phi_s))$ may be viewed as a sum
of the form (\ref{Wdef}) for a suitable $f$ which selects those
points with degree $k$. For each $x \in \BX$ we need to find coupled variables
$U_x$ and $V_x$ such that $U_x$ has the
distribution of $D_k (G(\Po_s ,\phi_s))$ and $1 + V_x$ has the
 conditional distribution of 
 $D_k (G(\Po_s \cup \{x\},\phi_s))  $ 
 given that $x$ has degree $k$,
 and such that $|U_x-V_x|$ is small (in probability). 
To do this we note that by the thinning theorem (see for example \cite{LP}),
 the point process of points of $\Po_s$ connected to $x$, and the point process
 of points of $\Po_s$ not connected to $x$, are independent.
To generate $V_x$ we need to condition the first point process
to have exactly $k$ points, which we can do by adding or
removing points from it, while leaving the second point process
unchanged. Since we modify only the first point process
(i.e., the points connected to $x$),
we may hope that the score $V_x$ obtained from the modified
(overall) point process is similar to the
 score $U_x$ obtained from the original
(overall) point process.
This is how we shall prove Theorem \ref{Th.k}.

For Theorem \ref{Th.cpts},
we view $N_k (G(\Po_s,\phi_s))$ as a sum
of the form (\ref{Wdef}) for a suitable $f$ which selects those
$k$-tuples of points forming a component.
 For each $\bx \subset \BX$ with $k$ elements,
 we need to find coupled variables
$U_\bx$ and $V_\bx$ such that $U_\bx$ has the
distribution of $N_k (G(\Po_s ,\phi_s))$ 
and $1+ V_\bx $ has the conditional distribution of 
 $N_k (G(\Po_s \cup \bx,\phi_s))$ 
given that the points of $\bx$ form a component. Again
 by the thinning theorem, the point process
of points of $\Po_s$ connected to $\bx$ is independent
of the rest of $\Po_s$, so to get $V_\bx$ we condition on
this point process having no elements, simply by removing
those points. Again this
is a small change, so the difference between
$U_\bx$ and $V_\bx$ again is small (in probability).

\section{Formal constructions of the IRG}
\label{secconstruct}
Let $s \in (0,\infty)$ and $n \in \N$.  Let $\phi \in \Phi$.
We now give a more formal definition of the graphs
$G(\Po_s,\phi)$ and $G(\cX_n,\phi)$,

We make the following particular choice of  mark space
$(\BM^*,\cM^*,\bm^*) $.
 Let $\Leb$ denote Lebesgue measure on $[0,1)$.
Let $\BM^* = [0,1)^{\N_0}$ with $\cM^*$ 
 the product Borel 
$\sigma$-algebra,  and with $\bm^* := \otimes_{n=0}^\infty \Leb$,
so that a random element of $\BM^*$ with distribution $\bm^*$ is
a sequence of independent uniform$(0,1)$ random variables
indexed by $\N_0$.

Now (and for the rest of this paper) taking 
 $ (\BM,\cM,\bm) = (\BM^*,\cM^*,\bm^*) $,
let $(X_i,T_i)_{i \geq 1}$, $(\PPrv_s)_{s >0}$ and $(\tau,\tau_1,\tau_2,\ldots)$
be as in the preceding section.  Given $n \in \N$ and $s >0$,
let the point processes $\eta_s$ and $\xi_n$ be as given by
  (\ref{etaxi}). 
Thus  $\eta_s$ is a Poisson point process in $\BX$ with mean measure $s \mu$ 
and with each point marked with a sequence of independent 
uniform$[0,1)$ variables
indexed  by the nonnegative integers. 
Similarly, $\xi_n$ is a  binomial point process in $\BX \times \BM^*$. 
We write $\bS^*$ for $\bS(\BX \times \BM^*)$.

To ease notation, we shall also assume from now on
that the probability measure $\mu$ on $\BX$ 
is {\em diffuse}.  This ensures
 that the values of $X_1,X_2,\ldots$ are almost
surely distinct. However, this assumption is for notational 
convenience only; even without it, one can make sense of our
results either by allowing the set $\{X_1,\ldots,X_n\}$
to have multiplicities, or by using the attached 
marks $T_i$ (which are almost surely distinct)
to distinguish between different points $X_i$.

Set $\cX_n := \{X_1,\ldots,X_n\}$ and
 $\Po_s := \{X_1,\ldots,X_{\PPrv_s}\}$, the canonical projections
of $\xi_n$ and $\eta_s$ respectively onto $\BX$.
Define the graph $G(\cX_n,\phi)$ to have vertex set $\cX_n$
and to have an edge between vertices $X_i$ and
$X_j$, for $i, j \in [n]: = \{1,2,\ldots,n\}$ with $i < j$,
 if and only if $T_{i,j} \leq \phi(X_i,X_j)$,
where $T_i = (T_{i,0},T_{i,1},T_{i,2} \ldots)$. 
Let $G(\Po_s,\phi)$ be the graph $G(\cX_{\PPrv_s},\phi)$.  

This is one way to formally define
the random graphs with the properties described
more informally in the Introduction. It has the
advantage that $G(\cX_n,\phi)$ is a subgraph of $G(\Po_s,\phi)$ whenever
$n \leq \PPrv_s$ and
$G(\Po_s,\phi)$ is a subgraph of $G(\cX_n,\phi)$ whenever
$\PPrv_s \leq n$,
 which is useful for coupling arguments.
However, it has the disadvantage that the edge-set of
the graph $G(\Po_s,\phi)$ is not invariant under permutation
of the order in which the marked points $(X_1,T_1),\ldots,(X_{\PPrv_s},T_{\PPrv_s})$
are listed. Therefore we define a further graph which
has the same distribution but also satisfies this permutation-invariance.
 This will be useful in applying
 Theorem \ref{thSteink} in the 
 proof of Theorems \ref{Th.k}, \ref{Th.cpts}
and \ref{thclump}.

Given $\xi \in \bS^*$ and $\phi \in \Phi$, define the graph $G_\phi(\xi)$ as
follows. 
If there exist distinct $(x,t_0,t_1\ldots) \in \xi$ and
$(y,u_0,u_1,\ldots) \in \xi$ with $t_0 = u_0$,
  then set $G_\phi(\xi)$ to be the empty graph.
Otherwise, we can write $\xi $ uniquely as
$$
\xi = \cup_{i=1}^{|\xi|}  \{(x_i,t_{i,0},t_{i,1},t_{i,2},\ldots)\} 
$$
with each $x_i \in \BX$ and
with $t_{1,0} < t_{2,0} < \cdots < t_{|\xi|,0}$. Let
$G_\phi(\xi)$ have vertex set $\{x_1,\ldots,x_{|\xi|}\}$, 
 and for each $i < j \leq |\xi|$ let $G_\phi(\xi)$  have an edge
connecting $x_i$  to $x_j$ if and only if
$t_{i,j}  \leq \phi(x_i,x_j)$.
In other words, we use the first coordinate
of the marks to determine the order in which we enumerate
the points of $\xi$; having done so, for $i <j$ we
use the $(j+1)$-st component of the mark attached
to the $i$-th point to decide whether to connect it
to the $j$-th point.

Let us say that two random  graphs $G$ and $G'$ have
the same distribution if any graph invariant of $G$ has
the same distribution as the same graph invariant evaluated
on $G'$. The following is immediate from the independence
of  the components of the marks $\tau_i$.
\begin{prop}
\label{propRCM}
For any finite  $\cX = \{x_1,\ldots,x_m\} \subset \BX$
the distribution of the random
 graph $G_\phi (\{(x_1,\tau_1),\ldots, (x_m,\tau_m)\})$ 
is the same as that of $G(\cX,\phi)$.
\end{prop}

In particular, the the distribution of $G_\phi(\eta_s)$ is the
same as that of $G(\Po_s,\phi)$, although they are not the same
graph because the set of edges is defined differently
for the two graphs. Likewise
$ G_\phi(\xi_n)$ has the same distribution as $G(\cX_n,\phi) $.

In the following lemmas we check measurability of
functions which will feature in the proof of Theorems 
\ref{Th.k} and \ref{Th.cpts}  respectively.
We use the following notation.  For $\ell \in \N$ and $i \in [\ell]$,
let $I_{\ell,i}$ denote the set of
$(x,t_0,t_1,t_2,\ldots) \in \BX \times [0,1)^{\N_0}$
such that  $(i-1)/\ell \leq t_0 < i/\ell$.

\begin{lemm}
\label{lemmeas1}
Let $j \in \N_0$.  Then the function 
$f: \BX \times \BM^* \times \bS^* \to \{0,1\}$ given by 
$$
f(x,\bt,\xi) := {\bf 1} \{ x {\rm ~has ~degree}~j {\rm ~in} 
~G_\phi(\xi \cup \{(x,\bt)\}) \}, ~~~~~~
 (x,\bt,\xi) \in \BX \times \BM^* \times \bS^* 
$$
 is measurable. 
\end{lemm}

{\em Proof.}
 For $\ell,n,m,i_1\ldots,i_n \in \N$
 with $ m \leq n \leq \ell$,
and  $i_1 < i_2 < \cdots < i_n \leq \ell$,
 let $A_{\ell,n,m,i_1,\ldots,i_n}$ denote the set of
$(x,\bt, \xi) \in \BX \times \BM^* \times \bS^*$ 
such that:
\benu
\item $(\xi \cup \{(x,\bt)\} ) (I_{\ell,i_k}) =1,$ for $ 1 \leq k \leq n$; 
\item $(\xi \cup \{(x,\bt)\} ) (I_{\ell,i}) = 0,$ for $ i \in [\ell]  \setminus
\{i_1,\ldots,i_n\}$;
\item $(x,\bt) \in I_{\ell,i_m}$;
\item
$\left( \sum_{k =1}^{m-1} \xi ( \{(y,\bu) \in I_{\ell,i_k} : u_m \leq \phi (x,y) \} ) \right) 
+ 
\sum_{k =m+1}^n \xi ( \{(y,\bu) \in I_{\ell,i_k} : 
t_k \leq \phi (x,y) \} ) 
=j,$ 
 where $\bu=(u_0,u_1,\ldots)$ and
 $\bt=(t_0,t_1,\ldots)$.
\eenu
Then each $  A_{\ell,n,m,i_1,\ldots,i_n}$ is measurable
in $\BX \times \BM^* \times \bS^*$, and
$$
f(x,\bt, \xi) = {\bf 1}\{(x,\bt,\xi) \in \cup_{  
 \ell,n,m,i_1,\ldots,i_n \in \N :
m \leq n \leq \ell, 
1\leq i_1 < i_2 < \cdots < i_n \leq \ell} 
A_{\ell,n,m,i_1,\ldots,i_n} \}
$$
which is a measurable function.
$\qed$

\begin{lemm}
\label{lemmeas2}
Let $k \in \N$.
Suppose the function $\tf: (\BX \times \BM^*)^k
 \times \bS^* \to \{0,1\}$ is given, for
 $(x_1,\bt_1,\ldots,x_k,\bt_k) \in (\BX \times \BM^*)^k$ and
$\xi \in \bS^* $,  
by 
\bean
\tf(x_1,\bt_1, \ldots,x_k,\bt_k,\xi) 
~~~~~~~~~~~~~~~~~~~~~~~~~
~~~~~~~~~~~~~~~~~~~~~~~~~
~~~~~~~~~~~
\\
 := {\bf 1} \left\{ \{x_1,\ldots,x_k\}  {\rm ~induces ~a ~component ~of ~}
G_\phi \left( \cup_{i=1}^k \{(x_i,\bt_i)\}
 \cup \xi \right) \right\}.
\eean
Then $\tf$ is measurable. 
\end{lemm}
{\em Proof.}
 For $\ell,n,m_1,\ldots,m_k,i_1 \ldots,i_n \in \N$
 with $\max_{1 \leq i \leq k} m_i \leq n \leq \ell$,
and $m_1,\ldots,m_k$ distinct, 
and  $1\leq i_1 < i_2 < \cdots < i_n \leq \ell$,
 let $A_{\ell,n,m_1,\ldots,m_k,i_1,\ldots,i_n}$ denote the set of
$(x_1,\bt_1,\ldots,x_k,\bt_k, \xi) \in (\BX \times \BM^*)^k \times \bS^*$
such that for some connected graph $\Gamma$ on vertex set $[k]$:
\benu
\item $\left( \cup_{h=1}^k \{(x_h,\bt_h)\} \cup \xi \right)
 (I_{\ell,i_j}) =1,$ for $ 1 \leq j \leq n$; 
\item $\left( \cup_{h=1}^k  \{(x_h,\bt_h)\} \cup \xi \right) 
(I_{\ell,i}) =0,$ for $ i \in [\ell]  \setminus
\{i_1,\ldots,i_n\}$;
\item $(x_h,\bt_h) \in I_{\ell,i_{m_h}}$ for $ 1 \leq h \leq k$;
\item
 for all $h \in [ k]$ and $ 1 \leq j < m_h$ 
with $j \notin \{m_1,\ldots,m_k\}$
 we have
$
\xi (\{ (y,\bu) \in I_{i_j}:   u_{m_h} >\phi(x_h,y)\} ) =1,
$ 
 where $\bu =(u_{0},u_{1},\ldots)$;
\item
for all $h \in [k]$ and $m_h < j \leq n$ 
with $j \notin \{m_1,\ldots,m_k\}$ we have
$ \xi (\{ (y,\bu) \in I_{i_j}:   t_{h,j} > \phi(x_h,y)\} ) =1$,
 where $\bt_h =(t_{h,0},t_{h,1},\ldots)$, and
\item
for all $h,h' \in [k]$ with $m_h < m_{h'}$
and $\{h,h'\}$ an edge of $\Gamma$,
we have $t_{h,m_h } \leq \phi (x_h, x_{h'})$.
\eenu
Then each $  A_{\ell,n,k,m_1,\ldots,m_k,i_1,\ldots,i_n}$ is measurable
in $(\BX \times \BM)^k \times \bS^*$, and
setting $[n]^k_{\neq}$ to be the set of
$(m_1,\ldots,m_k) \in [n]^k$ 
such that $m_1,\ldots,m_k$ are distinct, we have that
$\tf$
 is the indicator
of the set
\bean
 \cup_{  
\ell,n,i_1,\ldots,i_n \in \N,
(m_1,\ldots,m_k) \in [n]^k_{\neq}:
 n \leq \ell,
1\leq i_1 < i_2 < \cdots < i_n \leq \ell
} 
A_{\ell,n,m_1,\ldots,m_k,i_1,\ldots,i_n} 
\eean
which is  measurable.
$\qed$

\section{Proof of Theorem \ref{Th.k}}
\label{secproveThk}
\allco

For $s >0$, let $\eta_s$, $\Po_s$, $\PPrv_s$, $\tau,\tau_1,\tau_2,\ldots$
 be as in the preceding section.  For $\phi \in \Phi$ and  $i \in \N_0 $,
we can obtain from the Mecke equation (\ref{eqMecke}) that
\bea
\E \De_i(G_{\phi}(\eta_s)) = s \int_\BX  
\frac{ (s  \int \phi (x,y) \mu (dy))^i }{i!}
  \exp \left( -
 s \int \phi (x,y) \mu (dy)
 \right) \mu(dx).
\label{ENk}
\eea
Here is a brief explanation of (\ref{ENk}). In (\ref{eqMecke}), 
we take $k=1$ and use the space $\bS^*$ rather than $\bS$.
For $(x,\bt,\xi) \in \BX \times \BM^* \times \bS^*$,
our function $g((x,\bt),\xi)$ takes the value 1 if
$x$ has degree $i$ in the graph $G_\phi((x,\bt) \cup \xi)$,
and otherwise takes the value zero.  
Given $Z_s=n$, the point process $\Po_s$ has $n$ points, each of
which is independently connected to $x$ with probability
$\int \phi(x,y) \mu(y)dy$. Hence by the thinning property
of the Poisson distribution  \cite[Proposition 1.3]{LP},
the number of points of
$\Po_s$ connected to $x$ in $G_\phi ((x,\tau) \cup \eta_s)$ 
is Poisson  with mean $s \int \phi(x,y) \mu(dy)$.

In the sequel, other formulae for expectations
of numbers of vertices, or $k$-tuples of vertices,
having  certain properties in terms of the
graph $G_\phi(\eta_s)$ (or equivalently, the graph $G(\Po_s,\phi)$),
 will also be justified by the Mecke formula. These arguments also
justify (3.5) of \cite{Soft}.

Recall from (\ref{kappadef}) that for $\phi \in \Phi$ we set
$\kappa(\phi): = \sup_{x \in \BX} \int \phi(x,y)\mu(dy)$.

\begin{lemm}
\label{lemlog}
Suppose the assumptions of Theorem \ref{Th.k} hold.  Then 
(\ref{EDiDj}) holds. Also
\bea
s \kappa(\phi_s) = \Theta(\log s),
\label{Theta2cond}
\eea
and
 \bea
 \E \De_{j}(G(\Po_s,\phi_s))  \to  \alpha. 
\label{0513h}
\eea
\end{lemm}
\begin{proof}
For all $s >0$, since we assume $\phi_s \in \Phi_\eps$, 
by (\ref{ENk}) we have for $i \in \N$ that 
\bea
\frac{  \eps s \kappa(\phi_s) }{i} \leq
\frac{\E \De_{i}(G(\Po_s,\phi_s))  }{\E \De_{i-1}(G(\Po_s,\phi_s))  } \leq 
\frac{ s  \kappa(\phi_s) }{i}.
\label{0513e}
\eea
We are assuming for some fixed $j \in \N_0$ that
as $s \to \infty$ we have 
\bea
\E \De_{\leq j} (G(\Po_s,\phi_s)) \to \alpha \in (0,\infty).
\label{0513f}
\eea
Hence $\E \De_{0}(G(\Po_s,\phi_s)) $ is
bounded, but also 
 $\E \De_{0}(G(\Po_s,\phi_s))  \geq se^{-s\kappa(\phi_s)} $ by
(\ref{ENk}), so
$ s \kappa(\phi_s) \to \infty$ as $s \to \infty$.

By (\ref{0513f}), $\E \De_j(G(\Po_s,\phi_s))$ remains bounded, and since
also $s\kappa(\phi_s) \to \infty$, if $j \geq 1$ then using (\ref{0513e}) we
have $\E \De_{j-1}(G(\Po_s,\phi_s))  \to 0$, and repeating the argument
we also have 
$ \E \De_{i}(G(\Po_s,\phi_s))   \to 0, $ for $i =0,1,\ldots,j-1$.
Hence by (\ref{0513f}) we have
(\ref{0513h}).

Using (\ref{0513h}) and (\ref{ENk}),  the assumed $\eps$-homogeneity
of the $\phi_s$, and the fact that $s\kappa(\phi_s) \to \infty$,
it is straightforward to show that
(\ref{Theta2cond}) holds. 
Then, using (\ref{0513e}) repeatedly, and (\ref{0513h}), we obtain
(\ref{EDiDj}).
\end{proof}

{\em Proof of Theorem \ref{Th.k}.}
Assume the assumptions of  Theorem \ref{Th.k} apply. Then (\ref{EDiDj})
follows from Lemma \ref{lemlog}, and (\ref{EDiDj})  gives us
the first part of (\ref{eq0lim}).

Next we show that $ \De_{j}(G(\Po_s,\phi_s)) \tod {\Prv}_\alpha  $.
To carry out the strategy outlined at the end of Section
\ref{secstein},
we shall apply Theorem \ref{thSteink} to the case
$\phi = \phi_s$ of the function
$f$ considered in Lemma \ref{lemmeas1}.
That is, for $(x,\bt,\xi) \in \BX \times \BM^* \times \bS^*$ we set
$f_s(x,\bt,\xi)$ to be the indicator of the statement
that $x$ has degree $j$ in $G_{\phi_s}( \xi \cup \{(x,\bt)\})$.
Then
 $\De_j(G_{\phi_s}(\eta_s)) = F_s(\eta_s)$, where $F_s$ is the function
$F$ obtained by using $f \equiv f_s$ in (\ref{Wdef}) (with $k=1$).

Let $s >0$ and $x \in \BX$. 
 If $j \geq 1$ suppose we also have
an extra sequence $(Y,Y_1,\ldots,Y_j)$ of
independent identically distributed random elements of
 $\BX $ with $\Pr[ Y \in dy] = \phi_s(x,y) \mu(dy) / \int \phi_s(x,z) \mu(dz)$,
 independent of $((X_i,T_i))_{i \geq 1}$ and
$(\PPrv_s)_{s >0}$ and  $(\tau,\tau_1,\tau_2,\ldots)$. Let 
$
\G_s $ be the graph
$ G_{\phi_s}(\eta_s \cup \{(x,\tau),(Y_1,\tau_1),\ldots,(Y_j,\tau_j) \})$,
with added edges from $x $ to each of
$Y_1,\ldots,Y_j$
(if not already included).
Let $\Po_{s,x}$ be the set of points of
$\Po_s$ that are connected to $x$ in this graph,
and set $\Po_s^x := \Po_s \setminus \Po_{s,x}$.
 
Let $U_x$ denote the
number of vertices of degree $j$ in 
the subgraph of $\G_s$ induced by vertex set $\Po_s$.
By Proposition \ref{propRCM},
 this graph has the distribution of $G(\Po_s, \phi_s)$,
so $U_x$ has the distribution of $F_s(\eta_s)$.

Now consider the subgraph of $\G_s$ induced by $\{x\} \cup \Po_s$.
This has the distribution of $G(\Po_s \cup \{x\} , \phi_s)$,
and $\Po_{s,x}$ is the set of vertices in this graph
lying adjacent to $x$.  Conditioning on $x$ having degree $j$
 amounts to conditioning on $|\Po_{s,x}| =j$.  We define a coupled point
 process $\Po^*$ (a subset of $\Po_s \cup \{Y_1,\ldots,Y_j\}$) with  the
 distribution of $\Po_s$ conditioned on $x$ having degree $j$, as follows.

If  $|\Po_{s,x}|  >j$ then we select $|\Po_{s,x}| -j$ elements of
$\Po_{s,x}$ uniformly at random and discard them from $\Po_s$
to get a point process $\Po^*$. If
$|\Po_{s,x}| <j$ we set
 $\Po^* := \Po_{s} \cup \cY_x$, where we set
 $$
\cY_x := \{Y_1,\ldots, Y_{j- |\Po_{s,x}|} \}  .
$$
Let $\G_s^*$ denote the subgraph of $\G_s$
induced by $\Po^* \cup \{x\} $.
Then $x$ has degree $j$ in $\G^*_s$.
  Let $V_x$ be the number of vertices in $\G_s^*$ having
degree $j$, other than $x$. This has the conditional distribution of
$F_s(\{(x,\tau)\} \cup \eta_s) -1$ given that
 $f_s( \{(x,\tau)\},\eta_s ) =1$.
This is because $\Po_{s,x}$ and $\Po_s^x$ are
independent Poisson processes, and 
conditioning on 
 $f_s( \{(x,\tau)\},\eta_s ) =1$ amounts to
conditioning on  the first of these two Poisson processes
having $j$ points.

If $|\Po_{s,x}|  >j$ then
$|U_x - V_x| \leq U'_x + V'_x $, where we set
$U'_x$ to be the number  of $y \in \Po_{s,x}$ such that
$y$ has $j$ neighbours in $\Po_s$, and $V'_x$ 
to be the number of pairs $(y,z)$ 
with $y \in \Po_{s,x}$, $z \in \Po_s$,
such that $z \neq y$,  $z$ is connected to $y$ and
 $z$ has at most $j$ neighbours in
$\Po_s^x$.  By the Mecke formula,
and the assumption that $\phi_s \in \Phi_\eps$ for all $s$,
writing $\kappa_s := \kappa(\phi_s)$ and recalling that $s \kappa_s \to \infty$
by (\ref{Theta2cond}), we have
\bea
\E U'_x &
= & \int s \phi_s(x,y) \left( \frac{ (s \int \phi_s (y,z) \mu(dz) )^j}{j! } \right)
 \exp \left(- \int s \phi_s(y,z)  \mu(dz) \right) \mu(dy)
\nonumber \\
 & = & O( (s \kappa_s)^{j+1} ) \exp(- \Theta(s \kappa_s)) = o(1),
\label{EUpr}
\eea
uniformly over $x \in \BX$. Also, using that $\pmax_s \leq 1-\eps$
for all $s$, we have that
\bea
\E V'_x & = & \sum_{i=0}^j \int s \phi_s(x,y) \int s \phi_s(y,z) 
 \left( \frac{ (s \int \phi_s (z,w) (1- \phi_s(x,w)) \mu(dw) )^i}{i! } \right)
\nonumber \\
& &   \times \exp \left(- \int s \phi_s(z,w) (1- \phi_s(x,w))
  \mu(dw) \right) \mu(dz) \mu(dy)
\nonumber \\
 & = & O((s \kappa_s)^{j+2} 
) \times \exp( - \Theta(s \kappa_s)) \to 0,
\label{EWpr}
\eea
uniformly over $x \in \BX$.

If $|\Po_{s,x}|  <j$ then
$|U_x - V_x | \leq \tilde{U}_x + \tilde{V}_x $,
where we set $\tilde{U}_x$ to be the number of $y \in \cY_x$
having at most $j$ neighbours in $\Po_s$, and $\tilde{V}_x$
 is the number of pairs $(y,z)$ with $y \in \cY_x$, $z \in \Po_s$
connected to $y$, and $z $ having at most $j$ neighbours
in $\Po_s$. Then
\bea
\E \tilde{U}_x \leq j \sum_{i=0}^j \int
\left( s  \kappa_s \right)^i \exp(- \Theta(s  \kappa_s ))
\left( \frac{\phi_s(x,y) 
}{\int \phi_s(x,w) \mu(dw) } 
\right) 
\mu(dy)
\to 0,
\label{EtU}
\eea
uniformly over $x \in \BX$.  Also
\bea 
\E \tilde{V}_x & \leq
 & \sum_{i=0}^{j}
 \int s 
\left( s \int \phi_s(z,w) \mu(dw) \right)^i \exp
\left(-\int s \phi_s(z,w) \mu(dw) \right) 
\nonumber \\
& & \times j \int \left( \frac{ \phi_s(x,y) \phi_s(z,y)}{\int \phi_s(x,w) \mu(dw) } \right)
 \mu(dy)
 \mu(dz) 
\nonumber
\\
& =
& 
 O((s  \kappa_s)^{j+1} 
) \times \exp(- \Theta(s  \kappa_s))
\to 0,
\label{EtV}
\eea
uniformly over $x \in \BX$.

Combining the estimates
(\ref{EUpr}), (\ref{EWpr}), 
(\ref{EtU}) and  (\ref{EtV}),
 and using Theorem \ref{thSteink}, gives us
 the first part of (\ref{eqpolim}), namely
 $ \De_{j}(G(\Po_s,\phi_s)) \tod {\Prv}_\alpha  $.

Now suppose $k > j$. Set $\beta_{k,s} := \E[ D_k(G(\Po_s,\phi_s)) ]$.
The argument just given, with $j$ replaced by $k$, shows that
$d_{TV}( D_k(G(\Po_s,\phi_s)),Z_{\beta_{k,s}} ) \to 0 $ as $s \to \infty$.
Also $\beta_{k,s} \to \infty$ by (\ref{EDiDj}), so that
$ (Z_{\beta_{k,s}} -  \beta_{k,s} )/\sqrt{\beta_{k,s}}  \tod {\cal N} $.
Hence
$ (D_k(G(\Po_s,\phi_s)) - \beta_{k,s} )/\sqrt{\beta_{k,s}}  \tod {\cal N} $,
which is (\ref{DiN}).

It remains to prove the second parts of (\ref{eqpolim}) and  (\ref{eq0lim}).
For $n \in \N$,
let $s(n) = n- n^{3/4}$ and $t(n) := n+ n^{3/4}$.
By Chebyshev's inequality,
 with high probability (i.e. with probability
tending to 1) we have 
$\PPrv_{s(n)} \leq n \leq \PPrv_{t(n)}$ so that
$\Po_{s(n)} \subset \cX_n \subset \Po_{t(n)}$.
Moreover, when this happens, 
$G(\cX_n,\phi)$ is the subgraph of $G(\Po_{t(n)},\phi)$
induced by $\cX_n$, and
$G(\Po_{s(n)},\phi)$ is the subgraph of $G(\cX_{n},\phi)$
induced by $\Po_{s(n)}$.

By   (\ref{ENk}), for $i \in \N_0$,
$$
\left( \frac{s(n)}{n} \right)^{i+1} \leq
\frac{E D_i(G(\Po_{s(n)},\phi_n))}{E D_i(G(\Po_{n},\phi_n))}
\leq 
 \sup_{x \in \BX} 
\left( \exp \left( n^{3/4} \int \phi_n (x,y) \mu(dy) \right) \right)
$$
 and by (\ref{Theta2cond}), both the upper and the lower bound
tend to 1. Therefore by (\ref{EDiDj}) and (\ref{0513h}) we have 
as $n \to \infty$ that
\bea
& \E \De_{i}(G(\Po_{s(n)},\phi_n))   \to 0, &  ~~~~~~~ i =0,1,\ldots,j-1; 
\label{0625a}
\\
& \E \De_{j}(G(\Po_{s(n)},\phi_n))   \to   \alpha. &
\label{0625b}
\eea

Let $k \in \N_0$ with $k \leq j$.  If $\Po_{s(n)} \subset \cX_n$, then
$$
\De_{\leq k} (G(\cX_n,\phi_n)) - \De_{\leq k} (G(\Po_{s(n)},\phi_n))  
= S_n - R_n,
$$
where $S_n$ denotes the number of points
of $\cX_n \setminus \Po_{s(n)}$ with degree at most
$k$ in $G(\cX_n,\phi_n)$, and $R_n$ is
the number of points of $\Po_{s(n)}$ with
degree at most $k$ in $G(\Po_{s(n)},\phi_n)$
but with degree
 at least $k+1$ in $G(\cX_n,\phi_n)$.

Let $S'_n$ denote the number of points 
of $\Po_{t(n)} \setminus \Po_{s(n)}$ that are connected to
 at most
$k$ points of $\Po_{s(n)}$,
 and let $R'_n$ be
the number of points of $\Po_{s(n)}$ with
degree at most $k$ in $G(\Po_{s(n)},\phi_n)$
but with degree at least $k+1$ in $G(\Po_{t(n)},\phi_n)$.  If  
  $\Po_{s(n)}  \subset \cX_n  \subset \Po_{t(n)}$,  
then $S_n \leq S'_n$ and $R_n \leq R'_n$. Hence
by Markov's inequality,
\bean
\Pr[ \De_{\leq k} (G(\cX_n,\phi_n)) \neq \De_{\leq k} (G(\Po_{s(n)},\phi_n)) ]
\leq  \Pr[ 
  \{\PPrv_{s(n)}  \leq n  \leq \PPrv_{t(n)}\} ^c ]   
+ \E S'_n + \E R'_n.
\eean
Now
\bean
 \E S'_n & = & 2 n^{3/4}
\sum_{i=0}^k \int_\BX \frac{ (s(n) \int \phi_n(
x,y) \mu(dy) )^i}{i!}
\exp \left(-s(n) \int \phi_n(x,y) dy \right) \mu(dx)
\\
& = & \left( \frac{ 2 n^{3/4} }{s(n)}  \right)
 \E \De_{\leq k} (
G(\Po_{s(n)} , \phi_n)),
\eean
which tends to zero by 
(\ref{0625a}) and (\ref{0625b}).
Also
\bean
 \E R'_n & \leq & 
s(n) \sum_{i=0}^k \int  \frac{( s(n) \int \phi_n(x,y) \mu(dy) 
)^i  }{i!} ~~~~~~~~~~~~~~~~~~~~~~~ 
\\
& & \times \exp \left(- s(n) \int \phi_n(x,y) \mu(dy) \right) 
\left\{ 2 n^{3/4} \int \phi_n(x,y)
\mu(dy) \right\} \mu(dx)
\\
&  \leq  & 2 n^{3/4} a_n \E D_{\leq k} (\Po_{s(n)}, \phi_n) ,
\eean
which tends to zero by
(\ref{Theta2cond}),
 (\ref{0625a}) and (\ref{0625b}).

Therefore with high probability we have
$ \De_{\leq k} (G(\cX_n,\phi_n)) = \De_{\leq k} (G(\Po_{s(n)},\phi_n)) $.
This holds both for $k=j$, and 
for $k =j-1$.
 Hence using the first part of
(\ref{eqpolim}) and the first part of (\ref{eq0lim})
we obtain the
second part of (\ref{eqpolim}) and of (\ref{eq0lim}).
$\qed$ 

\section{Proof of Theorems \ref{Th.moments} and \ref{Th.cpts}}
\label{secGamma}
\allco
Given $\phi \in \Phi$, given $k, \ell \in \N$, and given
$\bx = (x_1\ldots,x_k) \in \BX^k$ and $\by = (y_1,\ldots,y_\ell ) \in \BX^\ell$,
 set
\bea
\phi(\bx,\by) 
:= 1- \prod_{i=1}^k \prod_{j=1}^\ell (1- \phi(x_i,y_j)).
\label{gphidef}
\eea
We also write
$ \phi(\{x_1,\ldots,x_k\}, \{y_1,\ldots,y_\ell\}) $ for
$\phi(\bx,\by)$ (allowing multiplicities in the sets $\{x_1,\ldots,x_k\}$
and $\{y_1,\ldots,y_\ell\}$); it is
 the probability that there is at least one edge in the random
graph $G(\{x_1,\ldots,x_k,y_1,\ldots,y_\ell\},\phi) $ connecting
one of the vertices $x_i$ to one of the vertices $y_j$.  
If $k=1$ we write
$ \phi(x_1,\{y_1,\ldots,y_\ell\}) $ for
$ \phi(\{x_1\},\{y_1,\ldots,y_\ell\}) $.
Also, let
$h_\phi(\bx) $ or $ h_{\phi}(x_1,\ldots,x_k)$
 denote the probability
that $G(\{x_1,\ldots,x_k\};\phi) $ is connected;
 more precisely, let $h_\phi(\bx) := 1 $ if $k=1$ and otherwise let
\bea
h_{\phi}(x_1,\ldots,x_k) := \sum_{\Gamma}
 \prod_{\{\{i,j\} : \{i,j\} \in E(\Gamma) \} } 
\phi(x_i,x_j)  
 \prod_{\{\{i,j\} : i,j \notin E(\Gamma) \} } 
(1 - \phi(x_i,x_j) ),
\label{hphidef}
 \eea
where the sum is over all connected graphs $\Gamma$ on vertex set
$\{1,\ldots,k\}$, and
$E(G)$ denotes the set of edges of a graph $G$. 
By the Mecke formula (\ref{eqMecke}),
and the equality in distribution of
$ G(\Po_s,\phi)$ and $G_{\phi}(\eta_s)$ as discussed in 
 Section \ref{secconstruct},
\bea
\E N_k(G(\Po_s,\phi)) 
= \E N_k(G_{\phi}(\eta_s)) 
~~~~~~~~~~~~~~~~~~~~~~~~~~~~~~~~~~~~~~
~~~~~~~~~~~~~~~~~~~~~~~~~
\nonumber \\
=  \frac{
 s^k}{k!}
\int_{\BX^k} 
 h_{\phi}(x_1,\ldots,x_k) 
\exp \left( - s \int \phi(z; \{x_1,\ldots,x_k\} ) \mu( dz ) \right)
 \mu^k ( d(x_1,\ldots,x_k)).
~~
\label{0625A}
\eea

Now fix $k \in \N$ and $\alpha \in (0,\infty)$. 
Assume throughout this section
that $\phi_s \in \Phi $ for $s >0$, and
(unless explicitly stated otherwise) that
\bea
\lim_{s \to \infty}
\E N_k(G(\Po_s,\phi_s)) = \alpha.
\label{ENklim}
\eea

{\em Proof of Theorem \ref{Th.moments} (a).}
Assume that $\pmax_s = o(1/\log s)$. 
We shall use the method of moments.

For $n,\ell \in \N$ we write
$(n)_\ell$ for the descending factorial $n(n-1)\cdots (n-\ell+1)$.
Then $( N_k(G(\Po_s,\phi_s)) )_{\ell}$
is the number of ordered $\ell$-tuples of distinct
$k$-components of $G(\Po_s,\phi_s) $.  This 
equals the sum over all ordered $k\ell$-tuples
$x_{1,1}, \ldots, x_{1,k}, \ldots, x_{\ell,1},\ldots,
x_{\ell,k}$ of distinct points of ${\Po_s}$, of
the indicator of the event that for each $i \leq \ell$ the subgraph of
$G(\Po_s,\phi_s)$ induced by $x_{i,1}\ldots,x_{i,k}$ is connected
and these vertices are not connected to any other vertices of 
$G(\Po_s,\phi_s)$, divided by $(k!)^{\ell}$.
Hence by the Mecke formula,
\bea
\E [ 
( N_k(G(\Po_s,\phi_s)) )_{\ell} ]
= \frac{
s^{k\ell} }{(k!)^\ell} \int_{\BX^{k\ell} }
\left( \prod_{i=1}^\ell h_{\phi_s} (x_{i,1},\ldots,x_{i,k} ) \right)  
 u_s(x_{1,1} \ldots,x_{\ell,k})
\nonumber \\
\times
\exp \left( -s \int  \phi_s(z,\{x_{1,1},\ldots,x_{\ell,k} \}) \mu(dz) \right) 
 \mu^{k \ell} (d(x_{1,1},\ldots ,x_{\ell,k})), 
\label{0625B}
\eea
where we set $u_s (\{x_{1,1},\ldots,x_{\ell,k}\}) $ 
to be the probability that the graph
$G(\{x_{1,1},\ldots,x_{\ell,k}\},\phi_s) $ 
has no edge between any $x_{i,j}$ and $x_{i',j'}$
such that $i \neq i'$, that is,
\bea
 u_s(x_{1,1} \ldots,x_{\ell,k}) :=
 \prod_{ (i_1,j_1), (i_2,j_2) \in [\ell] \times [k] : i_1 < i_2 } 
(1- \phi_s (x_{i_1,j_1},x_{i_2,j_2} )).
\label{undef}
\eea
By our condition on $\pmax_s$ the value of $u_s(x_{1,1},\ldots,x_{\ell,k})$
tends to 1, uniformly over $(x_{1,1},\ldots,x_{\ell,k})$.  
Also,  by the union bound 
$\phi_s(z,\{x_{1,1},\ldots,x_{\ell,k} \} )  $ is bounded by
\linebreak
$\sum_{i=1}^\ell \phi_s(z,\{x_{i,1},\ldots,x_{i,k} \} )$.  
Therefore by (\ref{0625B}),
 writing just $N_k$ for $N_k(G(\Po_s,\phi_s))$, we have
\bean
~~ \E [ ( N_k)_{\ell} ]
\geq (1+o(1))   \left(  \frac{
s^{k} }{k!} \int_{\BX^{k} }
h_{\phi_s} (x_{1},\ldots,x_{k} )  
\right.
~~~~~~~~~~~~~~
 ~~~~~~~~~~~~~~~~~~~~~~~~~~~~~~~~~
\\
 ~~~~~~~~~~~~~~
 \times
\left.
\exp \left( -s \int  \phi_s(z,\{x_{1},\ldots,x_{k} \}) \mu(dz) \right) 
 \mu^{k } ( d(x_{1},\ldots ,x_{k})) \right)^\ell. 
\eean
By (\ref{ENklim})  
 and (\ref{0625A}), this
lower bound  for 
$\E [ ( N_k)_{\ell} ] $ tends to $\alpha^\ell$ as $s \to \infty$.

By the  Bonferroni bound and the union bound 
we have for $(z,x_{1,1},\ldots,x_{\ell,k}) \in \BX^{1+ \ell k}$
that
\bea
\phi_s(z,\{x_{1,1},\ldots,x_{\ell,k} \} )  & \geq &
\sum_{i=1}^\ell
\phi_s(z,\{x_{i,1},\ldots,x_{i,k} \} )  
\left(1 - \sum_{j=i+1}^\ell 
\phi_s(z,\{x_{j,1},\ldots,x_{j,k} \} ) \right)  
\nonumber \\ &
\geq &
\sum_{i=1}^\ell
\phi_s(z,\{x_{i,1},\ldots,x_{i,k} \} )  
(1- k \ell \overline{\phi}_s). 
~~~~~~~~
\label{0625c}
\eea
Therefore by (\ref{0625B}),
\bean
\E [ (N_k)_\ell] 
\leq \frac{ s^{k \ell} }{(k!)^\ell}
\int_{\BX^{k\ell}}
\left( \prod_{i=1}^\ell h_{\phi_s} (x_{i,1},\ldots,x_{i,k} ) \right)  
~~~~~~~~~~~~ ~~~~~~~~~~~~
~~~~~~~~~~~~
 ~~~~~~~~~~~~
\\
 ~~~~~~~~~~~~
\times
\exp \left( -s \sum_{i=1}^\ell  
\int  \phi_s(z,\{x_{i,1},\ldots,x_{i,k} \}) 
(1 - k \ell \pmax_s)
\mu(dz) \right) 
 \mu^{k \ell} (d(x_{1,1},\ldots ,x_{\ell,k}) ) 
\\
= \left( 
\frac{ s^k  }{k!}
\int_{\BX^{k}}
h_{\phi_s} (x_{1},\ldots,x_{k} )  
\right.
~~~~~~~~~~~~ ~~~~~~~~~~~~
~~~~~~~~~~~~ ~~~~~~~~~~~~
~~~~~~~~~~~~ 
~~~ 
\\
\times
\left.
\exp \left( -s 
\int  \phi_s(z,\{x_{1},\ldots,x_{k} \}) 
(1 - k \ell \pmax_s)
\mu(dz) \right) 
 \mu^{k } ( d(x_{1},\ldots ,x_{k}) )
 \right)^\ell. 
\eean
Since $\mu$ is a probability measure
and $x^{1-k \ell \pmax_s}$ is a concave function on $x \geq 0$,
 we obtain by Jensen's inequality
and the fact that $h_{\phi_s}(\cdot)^{1/(1- k \ell \pmax_s)}
\leq h_{\phi_s}(\cdot)$
that
\bean
\E [ (N_k)_\ell] 
 \leq \frac{s^{k^2 \ell^2 \pmax_s}}{k!^{\ell^2 k \pmax_s}}
 \left( \frac{s^k}{k!} 
\int 
h_{\phi_s}(x_1,\ldots,x_k) \right.
~~~~~~~~~~~~~~~~~~~~~~~~~~~~~~~~~~~~~~~
~~~~~~~~~~
\\
~~~~~~~~~~
~~~~~~~~~~
 \times
\left.
\exp \left(- s \int \phi_s(z,\{x_1,\ldots,x_k\}) \mu(dz) 
\right) 
\mu^k(d(x_1,\ldots,x_k)) 
\right)^{\ell(1- k \ell \pmax_s)} 
\eean
which tends to $\alpha^\ell$ by (\ref{0625A}), (\ref{ENklim}) and
 our assumption on $\pmax_s$.
Thus $\E[(N_k)_{\ell}] \to \alpha^\ell$, so
 by the method of moments
(see e.g. Theorem 1.22 of \cite{Boll})
 the result (i) follows.$\qed$ \\

{\em Proof of Theorem \ref{Th.moments} (b).}
 Assume now that $\pmax_{s} = o(s^{-1/2})$.
Again we use the method of moments.
Write $n$ for $\lfloor s \rfloor$, and
 set $N'_k =  N_k(G(\cX_n,\phi_n))$.
Then by (\ref{undef}) and the union bound,
\bean
\E[ (N'_k)_\ell ]
 & = & \frac{
(n)_{k \ell}
}{k!^\ell} 
 \int \cdots \int 
\left( \prod_{i=1}^\ell h_{\phi_n} (x_{i,1},\ldots x_{i,k}) \right)
u_n (\{x_{i,1},\ldots,x_{\ell,k}\}) 
\\
&&  \times \left(1 - \int \phi_n(z; \{x_{1,1},\ldots,x_{\ell,k}\} )
 \mu(dz) \right)^{n-k\ell} 
\mu(dx_{1,1}) \cdots \mu(dx_{\ell, k}) \\
& \geq
& (1+o(1)) \left(
\frac{
  n^k}{k!}
 \right)^\ell 
 \int \cdots \int 
\left( \prod_{i=1}^\ell h_{\phi_n} (x_{i,1},\ldots x_{i,k}) \right)
\\
&&  \times \left( 1- \sum_{i=1}^{\ell}  \int
 \phi_n(z,\{x_{i,1},\ldots,x_{i,k}\}) \mu(dz) \right)^n  
\mu(dx_{1,1}) \cdots \mu(dx_{\ell,k}).
\eean
Using the bound $1-x \geq \exp(-x -x^2)$ for small positive $x$, we
have for large $n$ that
\bean
\E[ (N'_k)_\ell ] & \geq & (1+o(1)) \left(
\frac{
  n^k}{k!}
 \right)^\ell 
 \int \cdots \int 
\left( \prod_{i=1}^\ell h_{\phi_n} (x_{i,1},\ldots x_{i,k}) \right)
\\
&&  \times 
\exp \left( - n 
\left( \sum_{i=1}^\ell  \int \phi_n(z,\{x_{i,1},\ldots,x_{i,k}\})
 \mu(dz) \right) 
\right.
\\
&&
~~~~~~~~~~~ 
  - n \left.
 \left( \sum_{i=1}^\ell
\sum_{j=1}^k
\int \phi_n(z,x_{i,j}) \mu(dz) \right)^2  \right)
\mu(dx_{1,1}) \cdots \mu(dx_{\ell,k}),
\eean
so that
\bean
\E[ (N'_k)_\ell ] \geq 
 (1+o(1)) \left( \frac{n^k 
}{k!} 
 \int
 h_{\phi_n} (x_{1},\ldots x_{k})  \right.
~~~~~~~~~~~~~~~~~~~~~~~~
~~~~~~~~~~~~~~~~~~~~~~~~
\\
 \times
\left.
 \exp \left(- n  \int \phi_n(z,\{x_{1},\ldots,x_{k} \}) \mu(dz) \right)
\mu^k(d(x_{1}, \ldots ,dx_{k})) \right)^\ell   \times \exp( -nk^2 \ell^2
  \pmax_n^2 ) 
\eean
which tends to $\alpha^\ell$ by
 (\ref{0625A}), (\ref{ENklim}), and
 the assumption that
 $\pmax_n= o(n^{-1/2})$.

Conversely, by the bound $e^x \geq 1+x$ for all $x \in \R$,
we also have
\bean
\E [(N'_k)_\ell] \leq 
\left( \frac{n^{k} }{k!}  \right)^\ell 
\int \cdots \int
\left( \prod_{i=1}^\ell h_{\phi_n} (x_{i,1},\ldots x_{i,k}) \right)
~~~~~~~~~~~~~~~~~~~~~~~~
~~~~~~~~~
\\
\times
\exp \left( -  (n- \ell k) \int \phi_n(z;\{x_{1,1},\ldots,x_{\ell,k}
\})  \mu(dz) \right) \mu(dx_{1,1}) \cdots \mu(dx_{\ell,k})
\eean
and using the Bonferroni bound as in (\ref{0625c}),
we obtain that
\bean
\E [(N'_k)_\ell]
\leq e^{k^2 \ell^2 \pmax_n}
\left( \frac{n^{k} }{k!}  \right)^\ell 
\int \cdots \int
\left( \prod_{i=1}^\ell h_{\phi_n} (x_{i,1},\ldots x_{i,k}) \right)
~~~~~~~~~~~~~~~~~~~~
~~~~~~
~~~~~~
\\
~~~~~~
 \times
\exp \left( -  n \sum_{i=1}^\ell
 \int \phi_n(z;\{x_{i,1},\ldots,x_{i,k} \})
(1 - k \ell \pmax_n )
  \mu(dz)
\right)
\mu(dx_{1,1})
\cdots
\mu(dx_{\ell,k}).
\eean
By Jensen's inequality, since $x^{1- k \ell \pmax_n}$ is a concave function
on $x \geq 0$ and $h_{\phi_n}(\cdot) \leq 1$, we have
\bean
\E [(N'_k)_\ell] \leq e^{k^2 \ell^2 \pmax_n}
n^{k^2 \ell^2 \pmax_n} \left( \frac{1 }{k!}  \right)^\ell 
 \left( n^{k } \int \cdots \int h_{\phi_n} (x_{1},\ldots x_{k}) 
\right.
~~~~~~~~~~~~~~~~
~~~~~
\\
~~~~~
~~~~~
 \times
\left.
\exp \left( -  n 
 \int \phi_n(z;\{x_{1},\ldots,x_{k} \})
\mu(dz)
\right)
\mu(dx_{1})
\cdots
\mu(dx_{k})
\right)^{\ell(1 - k \ell \pmax_n) } .
\eean
Hence, since our assumption $\pmax_n =o(n^{-1/2}) $ implies
that $\pmax_n = o(1/(\log n))$, we have
$\limsup \E [(N'_k)_\ell] \leq \alpha^\ell$. Therefore
$\E [(N'_k)_\ell] \to \alpha^\ell$, and the method of moments
gives us part (ii). $\qed$ \\

\begin{lemm}
\label{lem17}
Let $k \in \N$,
 $\eps > 0$ and $\phi_s \in \Phi_\eps$
for all $s > 0$. 
Set
$ \kappa_s := \kappa(\phi_s)
= \sup_{x \in \BX} \int \phi_s(x,y) \mu(dy)$.
Then
\bea
\E N_k(G(\Po_s,\phi_s)) = \Theta(s^{k}  \kappa_s^{k-1} )
\exp(-\Theta(s  \kappa_s)). 
\label{0619a}
\eea
\end{lemm}
\begin{proof}
Our starting point is (\ref{0625A}). We first bound 
the exponent in the exponential factor.
Let $z,x_1,\ldots,x_k \in \BX$ and $s >0$. Then 
by the union bound,   
$$
\phi_s(z,x_1) \leq \phi_s(z,\{x_1,\ldots,x_k\}) \leq \sum_{i=1}^k \phi_s(z,x_i),
$$
so integrating over $z$ and using the assumed $\eps$-homogeneity, we have
\bea
\eps  \kappa_s \leq \int_\BX  
 \phi_s(z,\{x_1,\ldots,x_k\}) \mu(dz) \leq k  \kappa_s.
\label{170904a}
\eea
Using (\ref{0625A}), this already proves the result for $k=1$,
so from now on assume
 $k \geq 2$.

 If $x_{i-1}$ is connected to $x_i$ for $2 \leq i \leq k$,
then $G(\{x_1,\ldots,x_k\};\phi_s)$ is connected; hence
$h_{\phi_s}(x_1,\ldots, x_k) \geq \prod_{i=2}^k  \phi_s(x_{i-1},x_i)$,
so using (\ref{0625A}) and (\ref{170904a}) we obtain that
\bea
\E N_k(G(\Po_s,\phi_s)) & \geq & \frac{s^k}{k!} \exp(- k s  \kappa_s)
\int_\BX \cdots \int_\BX
 \prod_{i=2}^k  \phi_s(x_{i-1},x_i)
\mu(dx_k) \cdots  \mu (dx_1)
\nonumber \\
& \geq & \frac{\eps^{k-1}  \kappa_s^{k-1}  s^k}{k!} \exp (- k s  \kappa_s). 
\label{170904b}
\eea
For an upper bound, observe that if $G(\{x_1,\ldots,x_k\};\phi_s)$
is connected, then there exists a permutation $\sigma$ of $\{1,\ldots,k\}$
with $\sigma(1)=1$
such that for $2 \leq i \leq k$ this graph has an edge from
$x_{\sigma(i)}  $ to
$\{x_{\sigma(1)},\ldots, x_{\sigma(i-1)} \}  $. Therefore, setting
\bea
h_{\phi_s}^*(x_1,\ldots,x_k) := 
\prod_{i=2}^k \phi_s(x_i,\{x_1,\ldots,x_{i-1}\})
\leq
\prod_{i=2}^k \sum_{j=1}^{i-1} \phi_s(x_i,x_j),
\label{hstar}
\eea 
we have by the union bound that
$
h_{\phi_s}(x_1,\ldots,x_k) 
 \leq \sum_{\sigma} h_{\phi_s}^* (x_{\sigma(1)},\ldots,x_{\sigma(k)} )
$
where the sum is over all such permutations. 
Hence using (\ref{0625A}) and (\ref{170904a}), 
setting  $y_j = x_{\sigma(j)}$
and taking the integrals 
in  order $y_k,\ldots,y_1$
 for each permutation $\sigma$, we have
\bean
\E N_k(G(\Po_s,\phi_s)) \leq \frac{s^k (k-1)!}{k!} \exp(-\eps s  \kappa_s)
\int_\BX \cdots \int_\BX 
h_{\phi_s}^*( y_1,\ldots,y_k) \mu(dy_{k})\cdots \mu(dy_1),
\eean
and then using the inequality in (\ref{hstar}) we have
\bean
\E N_k(G(\Po_s,\phi_s)) \leq s^k (k-1)! \exp(-\eps s  \kappa_s)
 \kappa_s^{k-1}.
\eean
Combined with (\ref{170904b}) this gives us (\ref{0619a}). 
\end{proof}

{\em Proof of Theorem \ref{Th.cpts}.}
Assume there exists $\eps > 0$ such that $\phi_s \in \Phi_\eps$
and $\pmax_s \leq 1-\eps$ for all $s$.
Set $\kappa_s := \kappa(\phi_s)$.
If $s \kappa_s$ remains bounded away from zero and infinity, then by
(\ref{0619a}) we have that $ \E N_k(G(\Po_s,\phi_s))  \to \infty$,
contradicting (\ref{ENklim}).
Hence for any sequence
of values of $s$ tending to infinity, there is  a subsequence
such that either $s \kappa_s \to 0$
or $s \kappa_s \to \infty$ as $s \to \infty$ along the subsequence.

Consider first the case with $s \kappa_s \to 0$.
 In this case, by (\ref{0619a}) and (\ref{ENklim})
we have $k \geq 2$ and $s^{k} \kappa_s^{k-1} = \Theta(1)$.  
Recalling that $H_k(G)$ denotes the number of
 connected induced subgraphs  
of a graph $G$ of order $k$, we have
$$
\E H_{k+1}(G(\Po_s,\phi_s)) 
= O( s^{k+1}  \kappa_s^k) =o(1).
$$  
Since $0 \leq H_k (G(\Po_s,\phi_s))  - N_k (G(\Po_s,\phi_s))  
\leq (k+1) H_{k+1}(G(\Po_s,\phi_s))$, we have
$$
\E [ H_k (G(\Po_s,\phi_s))  -  N_k (G(\Po_s,\phi_s)) ]  
\to 0,
$$
so by (\ref{ENklim}) we have that $\E H_k (G(\Po_s,\phi_s))  \to \alpha$.
 Hence by Theorem \ref{thclump} (a) we have $H_k(G(\Po_s,\phi_s)) \tod
{\Prv}_\alpha$.  By Markov's inequality
$\Pr[H_{k}(G(\Po_s,\phi_s)) - N_{k}(G(\Po_s,\phi_s)) \geq 1 ] \to 0$,
so we also have $N_k(G(\Po_s,\phi_s)) \tod
{\Prv}_\alpha$, which is
the first part of (\ref{limGPo}).

Suppose now that $s  \kappa_s \to \infty$.
Then by (\ref{0619a})
 and (\ref{ENklim}), one may readily deduce that
\bea
s \kappa_s = \Theta(\log s).
\label{Thetacond}
\eea

We seek to apply Theorem \ref{thSteink}. 
For $s >0$, let $\eta_s$, $\Po_s$, $\PPrv_s$, $\tau,\tau_1,\tau_2,\ldots$
 be as in Section \ref{secconstruct}.
Let $\tf_s$ be the function $\tf$  considered in  Lemma
 \ref{lemmeas2}, using the connection function $ \phi \equiv \phi_s$. That is, 
let $ \tf_s(x_1,\bt_1, \ldots,x_k,\bt_k,\xi)
$ be the indicator of the statement that
 $\{x_1,\ldots,x_k\}$ induces a component of
$G_{\phi_s}\left(
 \cup_{i=1}^k \{(x_i,\bt_i)\}
 \cup
\xi
\right) $.
Then with $\tF_s$ denoting the function 
 $F$ obtained by taking $f \equiv \tf_s$ in
 the definition  (\ref{Wdef}), we have that
$\tF_s(\eta_s) = N_k(G_{\phi_s}(\eta_s))$ which
has the same distribution as $N_k(G(\Po_s,\phi_s))$. 

For
$\bx = (x_1,\ldots,x_k) \in \BX^k$,
set $\cX = \{x_1,\ldots,x_k\}$ (allowing multiplicities).
Define the graph 
\bea
\tGG_s :=  G_{\phi_s} (\eta_s \cup \{(x_1,\tau_1),\ldots,(x_k,\tau_k)\}). 
\label{cGdef}
\eea
Let  $\Po_{s,\bx}$ be the set of points
of $\Po_s$ connected 
to at least one point of $\cX$ in 
$\tGG_s$
and let
$  \Po^{\bx}_s = \Po_s \setminus \Po_{s,\bx}$. 

The subgraph of  
$\tGG_s$
 induced by vertex set $\Po_s$ has
the same distribution as $G(\Po_s,\phi_s)$, and we shall
 refer to this subgraph as $G'(\Po_s,\phi_s)$. Likewise, we refer
to the subgraph of  $\tGG_s$ induced by vertex set $\Po_s^{\bx}$ as
  $G'(\Po_s^{\bx},\phi_s)$, and we refer to the subgraph of 
$\tGG_s$ induced by vertex set $\cX$ as $G'(\cX,\phi_s)$.

Let $U_\bx = N_k (G'(\Po_s ,\phi_s))$.
This has the same distribution as $\tF_s(\eta_s)$.

Let $V_\bx = N_k(G'(\Po_s^{\bx} ,\phi_s))$.
We claim that  this has the  same distribution
as conditional distribution of $\tF_s( \cup_{i=1}^k \{(x_i,\tau_i) \}  
 \cup \eta_s)-1$ 
given that 
$\tf_s( \cup_{i=1}^k \{(x_i,\tau_i)\},  \eta_s ) = 1$.
This is because by the Marking Theorem for Poisson processes
(see e.g. \cite{Kingman}),
the point processes $\Po_{s,\bx}$ and $\Po_{s}^{\bx}$
are independent and the statement that 
$\tGG_s$
has $\cX$ as the vertex set of a component is equivalent
to the statement that  (i)
 $\Po_{s,\bx}$ has no points and (ii) the graph $G'(\cX,\phi_s)$ is connected,
which is independent of the outcome of $G'(\Po_s^{\bx},\phi_s)$.

Also, $U_\bx - V_\bx = U'_\bx - V'_\bx$, 
where $U'_{\bx} $ denotes the number of $k$-components
of $G'(\Po_s,\phi_s)$ with at least one vertex in $\Po_{s,\bx}$, and
 $V'_{\bx}$ is the number of $k$-components of
$G'(\Po_s^{\bx},\phi_s)$ with at least one neighbour in $\Po_{s,\bx}$.

By the Mecke formula
 \bean
\E U'_{\bx} = \frac{s^k}{k!} \int_{\BX^k} \phi_s(\bx,\by) h_{\phi_s}
(\by) \exp \left( - s \int \phi_s(z,\by) \mu(dz) \right) \mu^k(d \by),
\eean 
where $\phi_s(\bx,\by)$ is given by (\ref{gphidef}).
Since we can choose the elements $y_1\ldots,y_k$ of 
$\by$ in an order such that $y_1$ is connected to
$\cX$ and for each  $j \geq 2$, $y_j$ is connected
to $\{y_1,\ldots,y_{j-1}\}$, by a similar argument to
the proof of Lemma \ref{lem17} we have
\bea
\E U'_\bx = O(s^k  \kappa_s^k) \times \exp \left(- \Theta(s  \kappa_s) \right),  
\label{0717a}
\eea
which tends to zero (uniformly over $\bx$).

 Now $V'_{\bx}$ is bounded by the number of pairs $(y,\bz)$ with
$y \in \Po_{s,\bx}$ and $\bz = (z_1,\ldots,z_k )$ with
$\{z_1,\ldots,z_k\}$ inducing a $k$-component of $G'(\Po_s^{\bx},\phi_n)$
and $\bz$ connected to $y$.
Hence by the Mecke equation,
\bea
\E V'_{\bx} & \leq &
s^{k+1}\int_\BX \int_{\BX^k} \phi_s(y,\bx) \phi_s(y,\bz) 
(1- \phi_s(\bx,\bz))
h_{\phi_s}(\bz)
\nonumber  \\
& & \times \exp \left( - s \int \phi_s(w,\bz)
(1- \phi_s(w,\bx)) \mu(dw) \right)  
 \mu^k(d\bz) \mu(dy)
\nonumber \\
& = & O((s  \kappa_s)^{k+1}) \exp( - \Theta(s  \kappa_s)), 
\label{0717b}
\eea
which tends to zero, uniformly over $\bx$;
 here we have used the assumption that
$\pmax_s \leq 1-\eps$, so that $1 - \phi_s (w,\bx) \geq  \eps^k$,
 for all $s$, $w$, $\bx$. Therefore $\E[|U_\bx - V_\bx|] = 
\E[|U'_\bx - V'_\bx|]  \to 0$,  uniformly over $\bx \in \BX^k$.
 Then we can use Theorem \ref{thSteink} to
get the first part of (\ref{limGPo}). 

Before completing the proof of part (a) of Theorem 
\ref{Th.cpts}, we prove part (b), so now instead of (\ref{ENklim}) we assume
 $\alpha_s := \E N_{k}(G(\Po_s,\phi_s)) \to \infty $,
but $\alpha_s =o(s)$. Then by (\ref{0619a}), 
for every sequence of values of $s$ tending to infinity, there is  
a subsequence such that either $s  \kappa_s \to \infty$
 or $k \geq 2$ and $s  \kappa_s \to 0$ as $s \to \infty$ along the subsequence.

In both cases,  the estimates (\ref{0717a}) and (\ref{0717b}) hold 
so by Theorem \ref{thSteink}
we have $d_{TV}(\LL(N_k(G(\Po_s,\phi_s))),\LL({\Prv}_{\alpha_s})) \to 0$.
 Since also $({\Prv}_{\alpha_s} - \alpha_s )/\sqrt{\alpha_s} \tod {\cal N}$,
it follows that $(N_k(G(\Po_s,\phi_s))
 - \alpha_s)/\sqrt{\alpha_s} \tod {\cal N}$.
That is, for every sequence of values of $s$ tending to infinity,
there exists a subsequence such that
 $(N_k(G(\Po_s,\phi_s))
 - \alpha_s)/\sqrt{\alpha_s} \tod {\cal N}$ as $s \to \infty$ along
the subsequence.
Hence
 $(N_k(G(\Po_s,\phi_s))
 - \alpha_s)/\sqrt{\alpha_s} $ converges in
distribution to $ {\cal N}$ as $s \to \infty$.
This completes the proof of part (b).

Now we return to part (a), so we go back to assuming (\ref{ENklim}).
 As in the corresponding part of the proof of Theorem \ref{Th.k},
for $n \in \N$
 set $s(n) = n- n^{3/4}$ and
 $t(n) = n+ n^{3/4}$.
Then $\Po_{s(n)} \subset \cX_n \subset \Po_{t(n)}$
with high probability, and also the point process
$\Po_{t(n)} \setminus \Po_{s(n)}$ is a Poisson
point process with mean measure $2n^{3/4} \mu(\cdot)$,
independent of $\Po_{s(n)}$.  By (\ref{Thetacond}),
$ n^{3/4}  \int \phi_n(y,x) \mu( dy) = o(1), $
uniformly over $x \in \BX$, and therefore by (\ref{ENklim})
the sequence $(\phi_n)_{n \in \N}$ satisfies
\bea
 \frac{s(n)^k}{k!} \int_{\BX^k} \exp \left(
- s(n) \int_{\BX} \phi_n (z,\bx) \mu(dz) \right) 
h_{\phi_n}(\bx)
\mu^k(d\bx) \to \alpha.
\lbl{prminus}
\eea

For $n, \ell \in \N$ with $1 \leq \ell \leq k$,
let $A_{n,\ell}$ be the event that
 at least one collection of $\ell $  of the added vertices of $\cX_n \setminus
\Po_{s(n)}$  lies in a $k$-component of $G(\cX_n, \phi_n)$.  Let $B_n$
be the event that at least one of the added vertices of $\Po_{t(n)} \setminus
\Po_{s(n)}$ is connected to one of the $k$-components
 of  $G(\Po_{s(n)},\phi_n) $. 

If $A_{n,\ell}$ occurs and $\Po_{s(n)} \subset \cX_n \subset \Po_{t(n)}$, 
then there is at least one pair
$(\cX,\cY)$, such that
  $\cX \subset \Po_{s(n)}$ has $k- \ell$
 elements, and
  $\cY \subset \Po_{t(n)} \setminus \Po_{s(n)}$ has $\ell$ elements, and
$\cX \cup \cY$ induces a connected subgraph of $G(\Po_{t(n)},\phi_n)$,
and  there is no connection between any vertex of $\cX \cup \cY$
and any vertex of $\Po_{s(n) } \setminus \cX$ (however, we do allow
other connections between  vertices of $\cX \cup \cY$
and other vertices of $\Po_{t(n) } \setminus \Po_{s(n)}$).
By the Mecke equation, the expected number of such pairs
equals
\bean
\frac{( 2 n^{3/4})^\ell}{\ell!} \int
\frac{ (s(n))^{k-\ell} }{(k-\ell)!}  \int 
h_{\phi_n} (x_1,\ldots,x_{k- \ell},y_1,\ldots,y_{ \ell })
~~~~~~~
~~~~~~~
~~~~~~~
~~~~~~~
\\
\times \exp \left(- s(n) \int \phi_n(z,\{x_1,\ldots,x_{k- \ell},y_1,\ldots,y_{\ell}\}) 
\mu(dz) \right)
\\
\mu^{k-\ell}(d(x_1,\ldots,x_{k-\ell})) 
 \mu^\ell (d(y_1,\ldots,y_\ell)) ,
\eean
 which tends to zero by \eq{prminus}, so $\Pr[A_{n,\ell}] \to 0$.
Also, the expected number of $k$-components in $G(\Po_{s(n)},\phi_n)$
 which are connected to at least one vertex
of $\Po_{t(n)} \setminus \Po_{s(n)}$ is at most
$$
\frac{(s(n) )^k}{k!} \int \exp \left(- s(n)  \int
\phi_n (z,\bx) \mu(dz) \right) h_{\phi_n}(\bx) 
 2 n^{3/4} k a_n \mu^k(d\bx), 
$$
and by \eq{prminus} and (\ref{Thetacond}) this 
 tends to zero.  Hence $\Pr[B_n] \to 0$.  By the
first part of (\ref{limGPo}) we have for $k \in \N$ that 
$ N_k( G(\Po_{s(n)},\phi_{n}))  \tod {\Prv}_\alpha$.  Also 
$$
\Pr[ N_k(G(\cX_n,\phi_n))  \neq
N_k(G(\Po_{s(n)},\phi_n)) ] \leq \Pr[ \cup_{\ell=1}^k A_{n,\ell}] + 
\Pr[B_n] + \Pr[ \{\PPrv_{s(n)} \leq n \leq \PPrv_{t(n)}\}^c],
$$
 which tends to  0, and the second part of \eq{limGBi} follows.
 $\qed$

\section{Number of edges}
\label{secGn}
\allco
For any graph $G$, according to our earlier notation
 $H_2(G)$ denotes the number of edges of $G$. Let $\phi \in \Phi$. Then
by the Mecke formula (\ref{eqMecke}), 
\bea
\E H_2(G(\Po_s,\phi))
=
\E H_2(G_{\phi}(\eta_s))
 = \frac{1}{2} \int_\BX \int_\BX \phi(x,y) s^2 \mu(dx) \mu(dy).
\label{0713a}
\eea
\begin{theo}
\label{th.edge}
Suppose $\phi \in \Phi$.
Set $ \alpha:=  \E H_2(G(\Po_s,\phi))  .$ Then
\bea
d_{TV}( H_2(G(\Po_s,\phi)), {\Prv}_\alpha )
\leq  (1 \wedge  \alpha^{-1} ) 
\int  \left( \int \phi(x,y) s \mu(dy) \right)^2 s \mu(dx).
\label{Nedgest}
\eea 
\end{theo}
{\em Proof.}
We shall use Theorem \ref{thSteink}.  Given $x,y \in \BX$
with $\phi(x,y) > 0$,
set $\tGG_s = G_{\phi}(\eta_s \cup \{(x,\tau_1),(y,\tau_2) \} )$.
Let $G'(\Po_s,\phi)$ denote the subgraph of
$\tGG_s $ induced by the vertex set $\Po_s$,
which has the same distribution as
 $G(\Po_s,\phi)$.

Set $U_{x,y}= H_2(G'(\Po_s,\phi))$ and
 $V_{x,y}= H_2(\tGG_s ) - {\bf 1} \{ \{x,y\} \in E(\tGG_s) \}$, 
where $E(G)$ denotes the set of edges of a graph $G$. 
Then $ V_{x,y}$ has the conditional distribution
of $H_2(\tGG_s)- 1 $ given that $\{x,y\} \in E(\tGG_s) $.
Also
 $V_{x,y} \geq U_{x,y}$ and
$$
\E[ V_{x,y}- U_{x,y} ] = \int_{\BX} (\phi(x,z) + \phi(y,z))
s \mu(dz) =: w(x,y).
$$
Hence by Theorem \ref{thSteink} the left hand side of (\ref{Nedgest})
is bounded by the expression
\bean
(1/2)(1 \wedge \alpha^{-1})
\int \int w(x,y) \phi(x,y) s^2 \mu(dx) \mu(dy)
\\
=
(1 \wedge \alpha^{-1})
\int \left( \int \phi(x,y) s \mu(dy) \right)^2 s \mu(dx), 
\eean
as required. $\qed $\\

For example, consider the geometric setting with $\mu $
having bounded, almost everywhere continuous
  density $f$ with respect to Lebesgue measure
on $\R^d$. Suppose $\phi_s(x,y) =
\phi(r_s^{-1} (x-y))$ for some fixed integrable, symmetric and
almost everywhere continuous $\phi$, and
 $r_s >0$ satisfying $s^2 r_s^d \to \beta$ for some $\beta >0$. 
Set $ \alpha_s := \E H_2(G(\Po_s,\phi_s))$.
Then  by (\ref{0713a}),
$$
\alpha_s \to \frac{\beta}{2} \left( \int \phi(z) dz\right) \int f(x)^2 dx =: \alpha. 
$$
By (\ref{Nedgest}) we have
$ d_{TV} (\LL(H_2(G(\Po_s,\phi_s))), \LL ({\Prv}_{\alpha_s}) ) = O(s^3r_s^{2d}) = O(s^{-1})$.

In particular $H_2(G(\Po_s,\phi_s)) \tod {\Prv}_\alpha $. This could 
possibly also be proved
by deriving a Poisson limit for
$H_2(G(\cX_n,\phi_n)) $ by adapting the argument in
\cite[Theorem 3.4]{PBk} to the RCM, and Poissonizing.
However, the Poissonization would seem to introduce an
error of at least $s^{-1/2}$ in the total variation
distance, so the rate of convergence would not be as good. \\

{\em Proof of Theorem \ref{thclump}.}
Let $k \in \N$ with $k \geq 2$. 
Let $\eps >0$ and assume $\phi_s \in \Phi_\eps$ for all
$s >0$. Set $G_s := G(\Po_s,\phi_s)$.
Set $\kappa_s = \sup_{x \in \BX} \int \phi_s(x,y)\mu(dy)$.
 Then with $h_\phi(\cdot)$ defined at (\ref{hphidef}),
by (\ref{eqMecke}) we have 
\bea
\E H_k(G_s) = \frac{1}{k!}
\int s^k h_{\phi_s}( x_1,\ldots,x_k)
\mu^k (d(x_1,\ldots,x_k))
= \Theta(s^k \kappa_s^{k-1}),
\label{EHk}
\eea
where the second relation is obtained similarly
to the proof of Lemma \ref{lem17}, using the fact that
for each connected graph
$\Gamma$ on $\{1,\ldots,k\}$ we can integrate the variables
$x_1,\ldots,x_k$
in an order $(x_{\sigma(1)},\ldots,x_{\sigma(k)})$ such that
for $2 \leq i \leq k$
 each successive   $\sigma(i)$ 
is connected in $\Gamma$ to one or more of $\sigma(1),\ldots,\sigma(i-1)$. 
Then each successive integral gives another factor of $\Theta(\kappa_s)$.   

 Let $\bx = (x_1,\ldots,x_k) \in \BX^k$, with
$x_1,\ldots,x_k$ distinct and with $h_{\phi_s}(x_1,\ldots,x_k) >0$.
Let $\tGG_s$ be the graph defined by (\ref{cGdef}),
but now conditioned on the subgraph induced
by $\{x_1,\ldots,x_k\}$ being connected. 
Denote
 by $G'(\Po_s,\phi_s)$
 the  subgraph of $\tGG_s$ induced by
$\Po_s$.
Set $U_\bx = H_k(G'(\Po_s , \phi_s))$ 
and 
$ V_\bx = H_k(\tGG_s) - 1.
$
Then $U_\bx$ has the distribution
of $H_k(G_s)$ and $1+ V_\bx$ has the
conditional distribution of
$H_k(G(\Po_s \cup \{x_1,\ldots,x_k\},\phi_s))$
given that $\{x_1,\ldots,x_k\}$ induces
a connected subgraph of this graph.

Now $V_\bx \geq U_\bx$ and we assert that
\bea
\E [V_\bx - U_\bx] = O \left(
 \sum_{j=1}^{k-1} (s^j \kappa_s^j) 
\right),
\label{0714a}
\eea
uniformly over $\bx \in \BX^k$.
To see this, observe that
  $V_\bx - U_\bx $ is the number of pairs $(\cX,\Y)$
with $\cX$ a non-empty subset of $\{x_1,\ldots,x_k\}$ and
$\cY$ a non-empty subset of $\Po_s$, such that
the subgraph of $\tGG_s$ induced by
  vertex set $\cX \cup  \cY$ is a connected graph of
order $k$ (to ease notation
we ignore the issue of multiplicities in this notation).
But then, similarly to the proof of Lemma \ref{lem17},
 we can take the successive elements  $y_i$ of $\cY$
 in an order such that each
of them is connected to at least one existing vertex from 
$\cX \cup \{y_1,\ldots,y_{i-1}\}$. 
Then each successive integral gives another factor of $O(s\kappa_s)$.   

For part (a), assume $\E H_k(G_s)  \to \alpha \in (0,\infty)$.
Then by (\ref{EHk}) and (\ref{0714a}) we have
$\E |V_\bx - U_\bx| = o(1)$, uniformly over $\bx \in \BX^k$.
Hence by Theorem \ref{thSteink}, $H_k(G_s) \tod {\Prv}_\alpha$, 
which is part (i).

 For part (b), set $\alpha_s = \E H_k(G_s)$,
and assume $\alpha_s \to \infty$ but $\alpha_s = o(s)$.
Then by (\ref{EHk}),  the assumption   $\alpha_s = o(s)$
implies that (recalling $k \geq 2$) we have
 $s \kappa_s \to 0$.  Then by (\ref{0714a}) and
Theorem \ref{thSteink}, we have
 that $d_{TV}(\LL(H_k(G_s)), \LL({\Prv}_{\alpha_s}))
\to 0$. Since also $({\Prv}_{\alpha_s} - \alpha_s)/\sqrt{\alpha_s} \tod 
{\cal N}$, we therefore have
 $(H_k(G_s) - \alpha_s)/\sqrt{\alpha_s} \tod {\cal N}$.
$\qed$

\section{U-statistics of a Poisson process}
\label{secU}
\allco

Let $k \in \N$, and
let $\bS_k(\BX) := \{\xi \in \bS(\BX): |\xi| =k\}$.
   Let $h: \bS_k(\BX)  \to \{0,1\}$ be measurable.
For $\xi \in \bS(\BX)$ set
$$
F(\xi) = \sum_{\psi \subset \xi:|\psi| =k }
 h(\psi). 
$$
Let $s >0$, and let $\eta$ be a Poisson process on $\BX$ with mean measure
$\lambda := s \mu$.  We seek to apply Theorem \ref{thSteink} to $W:= F(\eta)$
for this class of choices of $F$, called $U$-statistics of the Poisson 
process $\eta$.
Assume $\mu$ is diffuse; once again,
 this assumption 
is for notational convenience only.

For $\bx = (x_1,\ldots,x_k) \in \BX^k$ with $h(\{x_1,\ldots,x_k\}) =1$, 
set $U_{\bx} = F(\eta)$, and
$$
V_\bx = F \left(
\cup_{i=1}^k \{x_i\}
\cup 
\eta 
\right) -1.
$$ 
Clearly $V_\bx$ and $U_\bx$ have the required distributional
properties in the statement of Theorem \ref{thSteink}. 
Also $V_\bx \geq U_\bx$, and
\bean
|V_{\bx} - U_{\bx} |= \sum_{J \subset [k]: 1 \leq |J| < k }
~
\sum_{\psi \subset \eta: | \psi | = k-|J| }
 h \left( 
\psi \cup \{ x_i: i\in J\}
 \right) ,
\eean
so using the Mecke formula and integrating over $\bx$ we obtain
that
\bean
\int 
(\E |V_{\bx} - U_{\bx} | ) h(\bx) \lambda^k(d\bx)
= \sum_{\ell= 1}^{k-1} \frac{1}{(k-\ell)!} \binom{k}{\ell}
\int  \int h(\{x_1,\ldots,x_\ell,y_1, \ldots, y_{k- \ell}\})
\\
\times \lambda^{k-\ell}(d(y_1,\ldots,y_{k-\ell})) 
h(\{x_1,\ldots,x_k\}) 
\lambda^k(d\bx)
\\
= \sum_{\ell= 1}^{k-1}
\frac{k!}{\ell! (k-\ell)!^2}
\int  \left( \int h(\{x_1,\ldots,x_\ell,y_1, \ldots, y_{k- \ell}\})
\lambda^{k-\ell}(d(y_1,\ldots,y_{k-\ell})) \right)^2
\\
\times
 \lambda^\ell(d(x_1,\ldots,x_\ell))
\\
=: \gamma(h,\lambda).
\eean
By Theorem \ref{thSteink}, if $\E F(\eta) = \alpha$, then
$d_{TV}(\LL(F(\eta)),\LL({\Prv}_\alpha)) \leq  (1 \wedge \alpha^{-1}) 
\gamma(h,\lambda)/k!$. Also
$d_{W}(\LL(F(\eta)),\LL({\Prv}_\alpha)) \leq 3 (1 \wedge \alpha^{-1/2}) 
\gamma(h,\lambda)/k!$. 
%
This bound  is comparable to the one obtained in 
Theorem 7.1 of \cite{DST}. Our bound has an extra factor
 $1 \wedge \alpha^{-1/2}$ in front, which may make it better
when $\alpha$ is large. Also, unlike \cite{DST} we do not
make any topological assumptions on the measurable
space $\BX$. As remarked just after the statement of 
Theorem \ref{thSteink}, it is possible to extend that
result to the case where the measure $\lambda$ is
$\sigma$-finite, and hence to extend the above 
argument likewise, but we do not go into details here. \\

{\bf Acknowledgements.}
We thank the anonymous referees for some valuable comments
on a previous version of this paper.

Part of this work was done while attending
the programme Theoretical Foundations for Statistical Network Analysis
at the Isaac Newton Institute for Mathematical Sciences, Cambridge 
 (EPSRC Grant Number EP/K032208/1).

{}

\end{document}